\documentstyle{amsppt}

\magnification=\magstep1
\baselineskip=14pt
\parindent=0pt
\parskip=14pt

\vsize=7.2in
\voffset=-.4in


\def \today {\ifcase \month \or January\or February\or  
  March\or April\or May\or June\or July\or August\or 
  September\or October\or November\or December\fi\   
  \number \day, \number \year} 
\def\now{ 
\countdef\hours=1
\countdef\minutes=2
\count0=\time
\divide\count0 by 60
\hours=\count0
\count3=\count0
\multiply\count3 by 60
\minutes=\time
\advance\minutes by -\count3
\number\hours:\number\minutes}

\overfullrule=0pt

\define\A{{\Bbb A}}
\define\C{{\Bbb C}}
\define\F{{\Bbb F}}

\define\R{{\Bbb R}}
\define\Q{{\Bbb Q}}

\define\Z{{\Bbb Z}}


\redefine\H{\frak H}
\define\a{\alpha}
\redefine\b{\beta}

\define\e{\epsilon}
\redefine\l{\lambda}
\redefine\o{\omega}
\define\ph{\varphi}

\define\s{\sigma}
\redefine\P{\Phi}
\predefine\Sec{\S}
\redefine\L{\Lambda}

\define\back{\backslash}

\define\lra{\longrightarrow}
\redefine\tt{\otimes}
\define\scr{\scriptstyle}
\define\liminv#1{\underset{\underset{#1}\to\leftarrow}\to\lim}
\define\limdir#1{\underset{\underset{#1}\to\rightarrow}\to\lim}

\define\isoarrow{\ {\overset{\sim}\to{\longrightarrow}}\ }

\define\nass{\noalign{\smallskip}}

\define\und#1{\underline{#1}}


\define\CH{\widehat{\roman{CH}}}


\define\db{\bar\partial}

\redefine\ord{\text{\rm ord}}
\define\Ei{\text{\rm Ei}}
\redefine\O{\Omega}
\predefine\oldvol{\vol}
\redefine\vol{\text{\rm vol}}

\define\CT#1{\operatornamewithlimits{CT}_{#1}}   

\define\pr{\text{\rm pr}}


%

\redefine\div{\text{\rm div}}

%
%

\define\Spec{\text{\rm Spec}}

\define\Sym{\text{\rm Sym}}
\define\tr{\text{\rm tr}}

\font\cute=cmitt10 at 12pt

\define\kay{{\text{\cute k}}}

\define\sig{\text{\rm sig}}

\define\Diff{\text{\rm Diff}}

\define\OB{\Cal O_B}
\define\End{\text{\rm End}}
\define\diag{\text{\rm diag}}
\redefine\Im{\text{\rm Im}}

\redefine\Re{\text{\rm Re}}

\define\MW{\text{\rm MW}}

\define\GL{\text{\rm GL}}
\define\SL{\text{\rm SL}}

\define\PGL{\text{\rm PGL}}
\define\Wald{\text{\rm Wald}}

\define\triv{{1\!\!1}}
\define\Spf{\roman{Spf}}
\define\inv{\roman{inv}}

\define\rags{\rangle_{\!\!\text{\lower 4pt\hbox{\rm GS}}}}

\define\Ver{\roman{Vert}}
\define\gs#1#2{\langle \,#1,#2\,\rangle}
\define\degh{\widehat{\deg}}
\define\divh{\widehat{\div}}
\define\Pich{\widehat{\roman{Pic}}}
\define\MWt{\widetilde{\MW}}
\define\hfb{\hfill\break}

\define\phih{\widehat{\theta}}

\define\Zh{\widehat{\Cal Z}}

\define\thh{\widehat{\theta}}
\define\res{\roman{res}}

\define\ohat{\widehat{\o}}
\define\Gt{\widetilde{G}}
\define\Pt{\widetilde{P}}
\define\It{\widetilde{I}}
\define\Pht{\widetilde{\P}}

\baselineskip=14pt

\define\ZZ{\Cal Z}
\define\ZZh{\ZZ^{\roman{hor}}}

\define\qeq{\ \overset{??}\to{=}\ }

\define\baruchmao{\bf1}
\define\boecherer{\bf2}
\define\borchinventI{\bf3}
\define\borchinventII{\bf4}
\define\borchduke{\bf5}
\define\borchdukeII{\bf6}
\define\bost{\bf7}
\define\bostUMD{\bf8}
\define\boutotcarayol{\bf9}
\define\bruinierI{\bf10}
\define\bruinierII{\bf11}
\define\brkuehn{\bf12}
\define\bkk{\bf13}  
\define\cohen{\bf14}
\define\faltings{\bf15}
\define\funkethesis{\bf16}
\define\funkecompo{\bf17}
\define\garrettII{\bf18}
\define\gelbartps{\bf19}
\define\gsihes{\bf20}
\define\grossmont{\bf21}
\define\grossMSRI{\bf22}
\define\grosskeating{\bf23}
\define\grosskohnenzagier{\bf24}
\define\grosskudla{\bf25}
\define\grosszagier{\bf26}
\define\harriskudlaII{\bf27}
\define\harriskudlaIII{\bf28}
\define\hirzebruchzagier{\bf29}
\define\howe{\bf30}
\define\howeps{\bf31}
\define\hriljac{\bf32}
\define\jacquet{\bf33}
\define\jacquetwald{\bf34}
\define\kitaoka{\bf35}
\define\kottwitz{\bf36}                                                   
\define\annals{\bf37}
\define\bourbaki{\bf38}
\define\Bints{\bf39}
\define\kudlaICM{\bf40}
\define\kudlaMSRI{\bf41}
\define\ariththeta{\bf42}
\define\kmI{\bf43}
\define\krannals{\bf44}
\define\krHB{\bf45}
\define\krinvent{\bf46}
\define\krsiegel{\bf47}
\define\tiny{\bf48}   
\define\kryII{\bf49}   
\define\kryIII{\bf50}  
\define\kryIV{\bf51}  
\define\doubling{\bf52}
\define\kuehncrelle{\bf53}   
\define\li{\bf54}
\define\maillotroessler{\bf55}
\define\mcgraw{\bf56}
\define\niwa{\bf57}
\define\odatsuzuki{\bf58}
\define\pswald{\bf59}
\define\psrallis{\bf60}
\define\rallisinnerprod{\bf61}
\define\saito{\bf62}
\define\siegel{\bf63}
\define\shimurahalf{\bf64}
\define\shimuraconf{\bf65}
\define\soulebook{\bf66}
\define\tunnell{\bf67}
\define\waldshimura{\bf68}
\define\waldfourier{\bf69}
\define\waldsurvey{\bf70}
\define\waldcentral{\bf71}
\define\waldspurgerLHD{\bf72}
\define\waldquaternion{\bf73}
\define\yangden{\bf74}
\define\yangiccm{\bf75}
\define\yangMSRI{\bf76}
\define\zagierII{\bf77}
\define\zagier{\bf78}  
\define\zhang{\bf79} 

\centerline{\bf Modular forms and arithmetic geometry}
\medskip
\centerline{Stephen S. Kudla\footnote{Partially 
supported by NSF grant DMS-0200292 and by a Max-Planck Research Prize 
from the Max-Planck Society and Alexander von Humboldt Stiftung. }}
\vskip .5in

The aim of these notes is to describe some examples of 
modular forms whose Fourier coefficients involve quantities 
from arithmetical algebraic geometry. 
At the moment, no general theory of such forms exists, but the 
examples suggest that they should be viewed as a kind of 
arithmetic analogue of theta series and that there should 
be an arithmetic Siegel--Weil formula relating suitable 
averages of them to special values of derivatives of Eisenstein 
series.  We will concentrate on the case for which 
the most complete picture 
is available, the case of generating series 
for cycles on the arithmetic surfaces associated to Shimura curves over $\Q$, 
expanding on the treatment in \cite{\kudlaICM}. A more speculative 
overview can be found in \cite{\kudlaMSRI}.

In section 1, we review the basic facts about the 
arithmetic surface $\Cal M$ associated to 
a Shimura curve over $\Q$. These arithmetic surfaces are moduli 
stacks over $\Spec(\Z)$ of pairs $(A,\iota)$ over a base $S$,
where $A$ is an abelian scheme of relative dimension $2$ and $\iota$ is
an action on $A$ of
a maximal order
$O_B$ in an indefinite quaternion algebra $B$ over $\Q$. 
In section 2, we recall
the definition of the  arithmetic Chow group $\CH^1(\Cal M)$, following Bost, \cite{\bost},
and we discuss the metrized Hodge line bundle $\hat\o$
and the conjectural value of $\langle\hat\o,\hat\o\rangle$, where 
$\langle\ ,\ \rangle$ is the height pairing on $\CH^1(\Cal M)$.  
In the next two sections, we describe divisors $\ZZ(t)$, 
$t\in\Z_{>0}$, on $\Cal M$. These are defined as the locus of  
$(A,\iota,x)$'s where $x$ is 
a special endomorphism (Definition~3.1) of $(A,\iota)$ 
with $x^2=-t$. Since such an $x$ gives an action on $(A,\iota)$ of the order
$\Z[\sqrt{-t}]$ in the imaginary quadratic field $\kay_t=\Q(\sqrt{-t})$, 
the cycles $\ZZ(t)$ can be viewed as analogues of the familiar CM points on modular 
curves. 
In section 3, the complex points and hence the 
horizontal components of $\ZZ(t)$ are determined. In section 4, 
the vertical components of $\ZZ(t)$ are determined using 
the p-adic uniformization of the fibers $\Cal M_p$ of bad reduction of 
$\Cal M$. In section 5, we construct Green functions $\Xi(t,v)$ for the divisors 
$\ZZ(t)$, depending on a parameter $v\in \R^\times_{>0}$. When $t<0$, 
the series defining $\Xi(t,v)$ becomes a smooth function on $\Cal M(\C)$. These
Green functions 
are used in section 6 to 
define classes $\Zh(t,v)\in \CH^1(\Cal M)$, for $t\in \Z$, $t\ne0$, 
and an additional class $\Zh(0,v)$ is defined using $\hat\o$. 
The main result of section 6 (Theorem~6.3) says that generating series 
$$\thh(\tau) = \sum_{t\in \Z}\Zh(t,v)\,q^t, \qquad \tau=u+iv, \qquad q = e(\tau),$$
is the $q$--expansion of a (nonholomorphic) modular form of weight $\frac32$, 
which we call an arithmetic theta function. The proof of this result 
is sketched in section 7. The main ingredients are (i) the fact that the 
height pairing of $\thh(\tau)$ with various classes in $\CH^1(\Cal M)$, 
e.g., $\hat\o$, can be shown to be modular, and (ii) the result of Borcherds, 
\cite{\borchduke},
which says that a similar generating series with coefficients in the usual Chow 
group of the generic fiber $\roman{CH}^1(\Cal M_\Q)$ is a modular form 
of weight $\frac32$.  In section 8, we use the arithmetic theta function 
to define an arithmetic theta lift
$$\thh: S_{\frac32} \lra \CH^1(\Cal M), \qquad f\mapsto \thh(f) =
\langle\,f,\thh\,\rangle_{\roman{Pet}},$$ 
from a certain space of modular forms of weight
$\frac32$ to  the arithmetic Chow group. 
This lift is an arithmetic analogue of the classical theta 
lift from modular forms of weight $\frac32$ to automorphic forms 
of weight $2$ for $\Gamma = O_B^\times$. According to the results of 
Waldspurger, reviewed in section 9, the nonvanishing of this 
classical lift is controlled by a combination of local obstructions
and, most importantly, the central value $L(1,F)$ of the standard 
Hecke L-function\footnote{Here we assume that $F$ is a newform, so that
$L(1,F)=L(\frac12,\pi)$, where $\pi$ is the corresponding cuspidal automorphic 
representation.}
of the cusp form $F$ of weight $2$ coming from $f$ via the Shimura lift. 
In section 10, we describe a doubling integral representation 
(Theorem~10.1) of the Hecke 
L-function, involving $f$ and an Eisenstein series $\Cal E(\tau,s,B)$ 
of weight $\frac32$ and genus $2$. At the central point $s=0$, 
$\Cal E(\tau,0,B)=0$.  In the case
in which the root number of the L-function is $-1$, we obtain a formula (Corollary~10.2)
$$
\gs{\,\Cal E'_2(\pmatrix
\tau_1&{}\\{}&-\bar\tau_2\endpmatrix,0;B)\,}{\,\overline{f(\tau_2)}\,}_{\text{Pet,
$\tau_2$}}
= f(\tau_1)\cdot C(0)\cdot  
L'(\frac12,\pi),
$$
for an explicit constant $C(0)$, whose vanishing is controlled by local obstructions.
Finally, in section 11, we state a conjectural identity (Conjecture~11.1)
$$\gs{\thh(\tau_1)}{\thh(\tau_2)} \qeq \Cal E_2'(\pmatrix
\tau_1&{}\\{}&-\bar\tau_2\endpmatrix,0;B)$$
relating the height pairing of the arithmetic theta function 
and the restriction to the diagonal of the 
derivative at $s=0$ of the weight $\frac32$ Eisenstein series. This 
identity is equivalent to a series of identities of 
Fourier coefficients, (11.1), 
$$\align
&\gs{\Zh(t_1,v_1)}{\Zh(t_2,v_2)}\cdot q_1^{t_1}q_2^{t_2}\\
\nass 
&\qquad\qquad{}=\sum_{\matrix\scr T \in \Sym_2(\Z)^\vee\\ \scr \diag(T) = (t_1,t_2)\endmatrix}
\Cal E'_{2,T}(\pmatrix \tau_1&{}\\{}&\tau_2\endpmatrix,0;B).
\endalign
$$
Here 
$$\Sym_2(\Z)^\vee = \{ \ T =\pmatrix t_1&m\\m&t_2\endpmatrix\mid t_1, t_2\in \Z, \ m\in
\frac12\Z\ \}$$
is the dual lattice of $\Sym_2(\Z)$ with respect to the trace pairing. 
 We sketch the proof of
these identities in the case where
$t_1t_2$ is not a square (Theorem~11.2). As a consequence, we prove Conjecture~11.1 up to 
a linear combination of theta series for quadratic forms in one variable (Corollary~11.3).
Assuming that $f$ is orthogonal to such theta series, we can 
substitute the height pairing $\gs{\Zh(t_1,v_1)}{\Zh(t_2,v_2)}$
for the derivative of the Eisenstein series in the doubling identity and obtain,  
$$\gs{\thh(\tau_1)}{\thh(f)} = f(\tau_1)\cdot C(0)\cdot  
L'(\frac12,\pi),
$$
in the case of root number $-1$. This yields the arithmetic inner product 
formula
$$\gs{\thh(f)}{\thh(f)} = \gs{f}{f}\cdot C(0)\cdot L'(\frac12,\pi),$$
analogous to the Rallis inner product formula for the classical theta lift. 
Some discussion of the relation of this result to the Gross-Kohnen-Zagier 
formula, \cite{\grosskohnenzagier}, is given at the end of section 11.

Most of the results described here are part of a long term collaboration 
with Michael Rapoport and Tonghai Yang. I would also like to thank 
J.-B. Bost, B. Gross, M. Harris, J. Kramer, and U. K\"uhn for their 
comments and suggestions.  

\medskip
\medskip 

\centerline{\bf Contents}
\smallskip

{\obeylines\parskip=0pt
\Sec1.\  Shimura curves and arithmetic surfaces
\Sec2. Arithmetic Chow groups
\Sec3. Special cycles, horizontal components
\Sec4. Special cycles, vertical components
\Sec5. Green functions
\Sec6. The arithmetic theta series
\Sec7. Modularity of the arithmetic theta series
\Sec8. The arithmetic theta lift
\Sec9. Theta dichotomy, Waldspurger's theory
\Sec10. The doubling integral
\Sec11. The arithmetic inner product formula
}

\subheading{\Sec1. Shimura curves and arithmetic surfaces}

Let $B$ be an indefinite quaternion algebra over $\Q$ and let $D(B)$ be the product 
of the primes $p$ for which $B_p= B\tt_\Q\Q_p$ is a division algebra. For the 
moment, we allow the case $B= M_2(\Q)$, where $D(B)=1$. 
The three dimensional $\Q$--vector space
$$V=\{x\in B\mid \tr(x)=0\}\tag1.1$$
is equipped with the quadratic form $Q(x) = \nu(x) = -x^2$. Here $\nu$ (resp. $\tr$) 
is the reduced norm (resp. trace) on $B$; in the 
case $B=M_2(\Q)$, this is the usual determinant (resp. trace). 
The bilinear form associated to
$Q$ is  given by $(x,y) = \tr(x y^\iota),$
where $x\mapsto x^\iota$ is the involution on $B$ given by $x^\iota = \tr(x) - x$. 
The action of $H= B^\times$ on $V$ by conjugation, $h: x\to hxh^{-1}$,
 preserves the quadratic
form  and induces an isomorphism
$$H \isoarrow \roman{GSpin}(V),\tag1.2$$
where $\roman{GSpin}(V)$ is the spinor similitude group of $V$. 
Since $B$ is indefinite, i.e., since $B_\R=B\tt_\Q\R\simeq M_2(\R)$, 
$V$ has signature $(1,2)$.
Let
$$D= \{ \ w\in V(\C)\ \mid\ (w,w)=0,\ (w,\bar w)<0\ \}/\,\C^\times \ 
\subset\ \Bbb P(V(\C)),\tag1.3$$
so that $D$ is an open subset of a quadric in $\Bbb P(V(\C))$. Then the 
group $H(\R)$ acts naturally on $D$, and,
if we fix an isomorphism $B_\R\simeq M_2(\R)$,  
then there is an identification
$$\C\setminus \R \isoarrow D, \qquad z\mapsto w(z):=\pmatrix z&-z^2\\1&-z\endpmatrix\ \mod
\C^\times,\tag1.4
$$  
which is equivariant for the action of $H(\R)\simeq\GL_2(\R)$ on $\C\setminus \R$ by fractional
linear transformations. 

Let $O_B$ be a maximal order in $B$ and let $\Gamma = O_B^\times$. In the case
$B=M_2(\Q)$, one may take $O_B=M_2(\Z)$, so that $\Gamma=\GL_2(\Z)$. 
Also let 
$$K = (\hat O_B)^\times \subset H(\A_f),\tag1.5$$
where $\hat O_B = O_B\tt_\Z\hat\Z$, for $\hat\Z = \liminv{N}\ \Z/N\Z$. 
Then the 
quotient 
$$M(\C) = H(\Q)\back \bigg(\ D\times H(\A_f)/K\ \bigg) \simeq \Gamma\back D,\tag1.6$$ 
which should be
viewed as an orbifold, is the set of complex points of a Shimura curve $M$, 
if $D(B)>1$, or of
the modular  curve (without its cusp), if $D(B)=1$. 
From now on, we assume that $D(B)>1$, 
although much of what follows can be carried over for $D(B)=1$ with only 
slight modifications. The key point is to interpret $M$ as a moduli space. 

Let $\Cal M$ be the moduli stack over $\Spec(\Z)$ for pairs $(A,\iota)$ where $A$ 
is an abelian scheme over a base $S$ with an action $\iota:O_B\rightarrow \End_S(A)$ 
satisfying the determinant condition, \cite{\kryII}, \cite{\kottwitz}, \cite{\boutotcarayol},
$$\det(\iota(b);\roman{Lie}(A)) = \nu(b).\tag1.7$$
Over $\C$, such an $(A,\iota)$ is an abelian surface with $O_B$ 
action. For example, for $z\in D\simeq \C\setminus\R$, the isomorphism
$$\l_z:B_\R\simeq M_2(\R) \isoarrow \C^2, \qquad b\mapsto 
b\cdot \pmatrix z\\1\endpmatrix = \pmatrix w_1\\ w_2\endpmatrix\tag1.8$$ 
determines a lattice $L_z=\l_z(O_B)\subset \C^2$. The complex torus $A_z = \C^2/L_z$
is an abelian variety 
with a natural $O_B$ action given by left multiplication, and hence 
defines an object $(A_z,\iota)\in \Cal
M(\C)$.   Two points in $D$ give the same lattice if and only if they are in the 
same $O_B^\times$--orbit, and, up to isomorphism, every 
$(A,\iota)$ over $\C$ arises in this way. 
Thus, the construction just described gives an isomorphism
$$[\,\Gamma\back D\,]\isoarrow\Cal M(\C)\tag1.9$$  
of orbifolds, and $\Cal M$ gives a model 
of $M(\C)$ over $\Spec(\Z)$ with generic fiber
$$M = \Cal
M\times_{\Spec(\Z)}\Spec(\Q),\tag1.10$$
the Shimura curve over $\Q$.   

Since we are assuming that $D(B)>1$, $\Cal M$ is proper of relative 
dimension $1$ over $\Spec(\Z)$ 
and smooth over $\Spec(\Z[\,D(B)^{-1}])$.  We will ignore the stack
aspect from now on and simply view $\Cal M$ as an arithmetic surface over $\Spec(\Z)$. 

The surface $\Cal M$ has bad reduction at primes $p\mid D(B)$ and this reduction can 
be described via p-adic uniformization \cite{\boutotcarayol}, \cite{\krinvent}.
Let $\hat\O_p$ be Drinfeld's p-adic upper half plane. It is a formal 
scheme over $\Z_p$ with a natural action of $\PGL_2(\Q_p)$. Let $W =W(\bar\Bbb F_p)$ 
be the Witt vectors of $\bar\Bbb F_p$ and let
$$\hat\O_W = \hat\O_p\times_{\Spf(\Z_p)}\Spf(W)\tag1.11$$
be the base change of $\hat\O_p$ to $W$.
Also, let $\hat\O_W^{\bullet} = \hat\O_W\times \Z$, and let $g\in\GL_2(\Q_p)$ 
act on $\hat\O_W^\bullet$ by
$$g:(z,i)\to (\,g(z),i+\ord_p(\det(g))\,).\tag1.12$$
Let $B^{(p)}$ be the definite quaternion algebra over $\Q$ with invariants
$$\inv_\ell(B^{(p)}) =\cases -\inv_\ell(B)&\text{if $\ell= p, \ \infty$,}\\
\nass
\phantom{-}\inv_\ell(B)&\text{otherwise.}
\endcases\tag1.13
$$
Let $H^{(p)} = (B^{(p)})^\times$ and 
$$V^{(p)} = \{ \ x\in B^{(p)} \mid \tr(x)=0\ \}.\tag1.14$$
For convenience, we will often write $B'=B^{(p)}$, $H'=H^{(p)}$ and $V'=V^{(p)}$ 
when $p$ has been fixed. Fix isomorphisms
$$H'(\Q_p) \simeq \GL_2(\Q_p), \qquad \text{and} \qquad H'(\A_f^p)\simeq H(\A_f^p).\tag1.15$$
Let $\widehat{\Cal M}_p$ be the base change to $W$ of the formal completion of 
$\Cal M$ along its fiber at $p$. Then, the Drinfeld--Cherednik Theorem gives an
isomorphism of formal schemes over $W$
$$\widehat{\Cal M}_p \isoarrow 
H'(\Q)\back \bigg(\ \hat\O_W^\bullet \times H(\A_f^p)/K^p\ \bigg) \simeq \Gamma'\back
\hat\O_W^\bullet,\tag1.16$$
where $K = K_pK^p$ and 
$\Gamma' = H'(\Q)\cap H'(\Q_p)K^p$.
The special fiber $\hat\O_p\times_{\Z_p}\F_p$ of 
$\hat\O_p$ is a union of projective lines $\Bbb P_{[\L]}$
indexed by the vertices $[\L]$ of the building $\Cal B$ of $\PGL_2(\Q_p)$. 
Here $[\L]$ is the homothety class of the $\Z_p$--lattice $\L$ in $\Q_p^2$. 
The crossing points of these lines are double points indexed by the 
edges of $\Cal B$,  
and the action of $\PGL_2(\Q_p)$ on components is compatible with 
its action on $\Cal B$. Thus, the 
dual graph of the special fiber of $\widehat{\Cal M}_p$ is isomorphic 
to 
$\Gamma'\back \Cal B^\bullet$,
where 
$\Cal B^\bullet = \Cal B\times\Z$.

\subheading{\Sec2. Arithmetic Chow groups}

The modular forms of interest in these notes will take values in the
arithmetic Chow groups  of $\Cal M$. We will use the version of these
groups with real coefficients
defined  by Bost, \cite{\bost}, section~5.5. Let $\widehat Z^1(\Cal M)$ 
be the real vector space spanned by pairs $(Z,g)$, where 
$Z$ is a real linear combination of Weil divisors on $\Cal M$ 
and $g$ is a Green function for $Z$. In particular, if $Z$ is a Weil divisor, 
$g$ is a 
$C^\infty$ function on $\Cal M(\C) \setminus Z(\C)$, with a logarithmic
singularity along the  $Z(\C)$, and satisfies the Green equation
$$dd^c g+\delta_Z = [\o_Z],\tag2.1$$
where $\o_Z$ is a smooth $(1,1)$-form on $\Cal M(\C)$, and 
$[\o_Z]$ is the corresponding current. If $Z = \sum_i c_i\, Z_i$ is a real linear combination
of  Weil divisors, then $g=\sum_i c_i\,g_i$ 
is a real linear combination of such Green  functions. By construction, $\a \cdot
(Z,g) = (\a Z,\a g)$  for $\a\in \R$. The first arithmetic Chow group, 
with real coefficients, $\CH^1_\R(\Cal M)$, is then 
the quotient of
$\widehat Z^1(\Cal M)$ by the subspace spanned by pairs $\divh(f) = (\div(f),-\log|f|^2)$ 
where  
$f$ is a rational function on $\Cal M$, and $\div(f)$ is its divisor. 
Finally, we let
$$\CH^1(\Cal M) = \CH^1_\R(\Cal M)\tt_\R\C.\tag2.2$$
Note that restriction to the generic fiber yields a degree map
$$\deg_\Q:\CH^1(\Cal M) \lra \roman{CH}^1(M)\tt\C \overset{\deg}\to{\lra}\ \C.\tag2.3$$
The group $\CH^2(\Cal M)$ is defined analogously, and the arithmetic 
degree map yields an isomorphism
$$\degh:\CH^2(\Cal M)\isoarrow \C.\tag2.4$$
Moreover, there is a symmetric $\R$--bilinear height pairing\footnote{Here the symmetry 
must still be checked in the case of a stack, cf. section 4 of \cite{\kryII}.}
$$\langle\ \,,\ \rangle:\ \CH^1_\R(\Cal M)\times\CH^1_\R(\Cal M) \lra \R.\tag2.5$$
According to the index Theorem, cf. \cite{\bost}, Theorem~5.5, this pairing is nondegenerate
and  has signature $(+,-,-, \dots)$. We extend it to an Hermitian pairing  
on $\CH^1(\Cal M)$, conjugate linear in the second argument. 

Let $\Cal A$ be the universal abelian scheme over $\Cal M$ 
with zero section $\e$, and let
$$\o = \e^*\O^2_{\Cal A/\Cal M}\tag2.6$$
be the Hodge line bundle on $\Cal M$. 
We define the natural metric on $\o$ by letting
$$||s_z||^2_{\roman{nat}} = \left|\,\left(\frac{i}{2\pi}\right)^2
\int_{A_z} s_z\wedge\bar{s}_z\ \right|.\tag2.7$$
for any section $s:z\mapsto s_z$, where, for $z\in \Cal M(\C)$, 
$A_z$ is the associated abelian variety. 
As in section 3 of \cite{\kryII}, we set 
$$||\ ||^2 = e^{-2C}\,||\ ||^2_{\roman{nat}}\tag2.8$$
where $2C = \log(4\pi) + \gamma$, where $\gamma$ is Euler's constant.
The reason for this choice of normalization is explained in the 
introduction to \cite{\kryII}. 
The 
pair $\ohat = (\o,||\ ||)$ defines an element of $\Pich(\Cal M)$, the group of 
metrized line bundles on $\Cal M$. We write $\ohat$ for the image of this class 
in $\CH^1_\R(\Cal M)$ under the natural map $\Pich(\Cal M)\rightarrow 
\CH^1_\R(\Cal M)$. 

The pullback to $D$ of the restriction of $\o$ to $\Cal M(\C)$ 
is trivialized by the section $\a$
defined as follows. For $z\in D$, let $\a_z$ be the holomorphic 
$2$--form  
on $A_z = \C^2/L_z$ given by 
$$\a_z = D(B)^{-1} (2\pi i)^2\, dw_1\wedge dw_2\tag2.9$$
where $w_1$ and $w_2$ are the
coordinates on the  right side of (1.6). Then 
$$\o_\C = \big[\ \Gamma\back (D\times \C)\,\big],\tag2.10$$
where the action of $\gamma\in \Gamma$ is given by
$$\gamma: (z,\zeta) \longmapsto (\gamma(z), (cz+d)^2\,\zeta).\tag2.11$$
Thus, on $\Cal M(\C)$, $\o$ is isomorphic to the canonical bundle $\O^1_{\Cal M(\C)}$, 
under the map which sends $\a_z$ to $dz$. The resulting metric on $\O^1_{\Cal M(\C)}$ is
$$\align
||dz||^2 &= ||\a_z||^2 = e^{-2C}\, \left|\,\left(\frac{i}{2\pi}\right)^2
\int_{A_z} \a_z\wedge\bar{\a}_z\ \right|\\
\nass
\nass
{}&= e^{-2C}\,(2\pi)^{-2}\,(2\pi)^4 \,D(B)^{-2}\,\vol(M_2(\R)/O_B)\cdot\Im(z)^2\tag2.12\\
\nass
\nass
{}&= e^{-2C}\,(2\pi)^2\cdot \Im(z)^2.\\
\endalign
$$

In \cite{\kryII}, it was conjectured that
$$\gs{\ohat}{\ohat} \overset{??}\to{=}\ \zeta_{D(B)}(-1)\bigg[\ 2\frac{\zeta'(-1)}{\zeta(-1)} +1 
-2C -\sum_{p\mid D(B)} \frac{p\log(p)}{p-1}\ \bigg],\tag2.13$$  
where 
$$\zeta_{D(B)}(s) = \zeta(s)\cdot \prod_{p\mid D(B)} (1-p^{-s}),\tag2.14$$
and $2C=\log(4\pi)+\gamma$, as before. In the case $D(B)=1$, i.e., for a modular 
curve, the analogous value for $\gs{\ohat}{\ohat}$ was established, independently, 
by Bost, \cite{\bostUMD}, and K\"uhn, \cite{\kuehncrelle}, cf. the introduction 
to \cite{\kryII} for a further discussion of normalizations. It will be convenient 
to define a constant $\bold c$ by
$$\frac12\,\deg_\Q(\ohat) \cdot \bold c := \gs{\ohat}{\ohat} - \text{RHS of (2.13)}.\tag2.15$$
In particular, $\gs{\ohat}{\ohat}$ has the conjectured value if and only if $\bold c=0$. 
It seems likely that this conjecture can be proved using recent work of 
Bruinier, Burgos and K\"uhn, \cite{\brkuehn}, on heights of curves on Hilbert modular surfaces. 
Their work uses an extended theory of arithmetic Chow groups, 
developed by Burgos, Kramer and K\"uhn, \cite{\bkk}, which allows metrics with
singularities of the type which arise on compactified Shimura varieties. 
In addition, they utilize results of Bruinier,
\cite{\bruinierI}, 
\cite{\bruinierII}, concerning Borcherds forms. Very general conjectures about 
such arithmetic degrees and their connections with the equivariant arithmetic Riemann--Roch 
formula have been given by Maillot and Roessler, \cite{\maillotroessler}. 
Some connections of arithmetic degrees with Fourier coefficients of derivatives 
of Eisenstein series are discussed in \cite{\Bints}. 

\subheading{\Sec3. Special cycles, horizontal components}

Just as the CM points on the modular curves are constructed as the 
points where the corresponding elliptic curves have additional endomorphisms, 
cycles on $\Cal M$ can be defined by imposing additional endomorphisms as follows. 
\proclaim{Definition 3.1} For a given $(A,\iota)$, the space of 
{\bf special endomorphisms} of $(A,\iota)$ is
$$V(A,\iota) = \{\ x\in \End_S(A)\ \mid\ \tr(x)=0 \text{ and } 
x\circ\iota(b) = \iota(b)\circ x, \ \forall b\in
O_B\
\}.\tag3.1$$
Over a connected base, this free $\Z$--module of finite rank 
has a $\Z$--valued quadratic form $Q$ given by
$$-x^2 = Q(x)\cdot \roman{id}_A.\tag3.2$$
\endproclaim

\proclaim{Definition 3.2} 
For a positive integer $t$, let $\ZZ(t)$ be the moduli stack over $\Spec(\Z)$ of
triples $(A,\iota,x)$, 
where $(A,\iota)$ is as before and $x\in V(A,\iota)$ 
is a special endomorphism with $Q(x) = t$.
\endproclaim
There is a natural morphism 
$$\ZZ(t) \lra \Cal M, \qquad (A,\iota,x) \mapsto (A,\iota),\tag3.3$$
which is unramified, and, by a slight abuse of notation, we write $\ZZ(t)$ for the 
divisor on $\Cal M$ determined by this morphism.  

Over $\C$, such a triple $(A,\iota,x)$ is an abelian surface with an action 
of $O_B\tt_\Z \Z[\sqrt{-t}],$ i.e., with additional `complex multiplication' 
by the order $\Z[\sqrt{-t}]$ in the imaginary quadratic field $\kay_t=\Q(\sqrt{-t})$. 
Suppose that $A \simeq A_z = \C^2/L_z$, so that, by (1.8), the tangent space $T_e(A_z)$ 
is given as $B_\R\simeq
\C^2 = T_e(A_z)$.  Since the lift $\tilde x_\R$ of $x$ in the diagram
$$\matrix B_\R\hskip -36pt &\simeq&\C^2&\lra&A_z\\
\nass
\tilde x_\R=r(j_x)\downarrow\phantom{\tilde x_\R=r(j_x)}\hskip -36pt &{}&\tilde x\downarrow
\phantom{\tilde x}&{}&x\downarrow\phantom{x}\\
\nass
B_\R\hskip -36pt  &\simeq&\C^2&\lra&A_z\endmatrix\tag3.4
$$
commutes with the left action of $O_B$ and carries $O_B$ into itself, it is given by right
multiplication $r(j_x)$ by an element $j_x\in O_B\cap V$ with $\nu(j_x) = t$. 
Since the map $\tilde x$ is holomorphic, it follows that $z\in D_x$, the fixed 
point set of $j_x$ on $D$. To simplify notation, we will 
write $x$ in place of $j_x$. Then, we find that
$$\ZZ(t)(\C) = \big[\ \Gamma\back D_t\ \big],\tag3.5$$
where
$$D_t = \coprod_{\matrix\scr x\in O_B\cap V\\ \scr Q(x)=t\endmatrix} D_x.\tag3.6$$
In particular, $\Cal Z(t)(\C)$, which can be viewed as a set of CM points 
on the Shimura curve $M(\C)$, is nonempty if and only if the
imaginary quadratic field $\kay_t$ embeds in $B$.

The horizontal part $\ZZh(t)$ of $\ZZ(t)$ is obtained by taking the 
closure in $\Cal M$ of these CM--points, i.e., 
$$\ZZh(t) := \overline{\ZZ(t)_\Q}.\tag3.7$$

Finally, if we write $4t=n^2d$ where $-d$ is the discriminant of the 
imaginary quadratic field $\kay_t=\Q(\sqrt{-t})$, then the degree of the divisor
$\ZZ(t)_\Q$ is given by
$$\deg_\Q\ZZ(t) = 2 \delta(d,D(B))\,H_0(t,D(B)),\tag3.8$$
where
$$\delta(d,D(B)) = \prod_{p\mid D(B)} (1-\chi_d(p)),\tag3.9$$
and 
$$\align
H_0(t,D(B)) &= \sum_{c|n} \frac{h(c^2d)}{w(c^2d)}\tag3.10\\
\nass
{}&=\frac{h(d)}{w(d)}\cdot\bigg(
\sum_{\matrix\scr c\mid n\\ \scr (c,D(B))=1\endmatrix} c\,\prod_{\ell\mid
c}(1-\chi_d(\ell)\ell^{-1})\,\bigg).
\endalign
$$
Here $h(c^2d)$ is the class number of the order $O_{c^2d}$ in $\kay_t$ 
of conductor $c$, $w(c^2d)$ is the
number  of roots of unity in $O_{c^2d}$, 
and $\chi_d$ is the Dirichlet character for the field $\kay_t$. 
Note that, in this formula, we are counting points on the orbifold $[\Gamma\back D]$, 
so that each point $\pr(z)$, $z\in D$, is counted with multiplicity $e_z^{-1}$ where 
$e_z = |\Gamma_z|$ is the order of the stabilizer of $z$ in $\Gamma$. 
For example, suppose that $z\in D_x$ for $x\in V\cap O_B$ with $Q(x)=t$. Then, 
since $\Z[x]\simeq \Z[\sqrt{-t}]$ is an order of conductor $n$, 
$\Q[x]\cap O_B \supset \Z[x]$ is an order of conductor $c$ for some $c\mid n$, 
and $e_z = w(c^2d) = |\ (\Q[x]\cap O_B)^\times|$.

\subheading{\Sec4. Special cycles, vertical components}

In this section, we describe the vertical components of the 
special cycles $\ZZ(t)$ in some detail, following \cite{\krinvent}.
In the end, we obtain a `p-adic uniformization', (4.27), quite analogous to 
the expression (3.5) for $\ZZ(t)(\C)$.  
We first review the construction of 
the p-adic uniformization isomorphism (1.16). 

Let $\Bbb B$ be the division quaternion algebra over $\Q_p$ and let $\OB$ 
be its maximal order. Let $\Z_{p^2}$ be the ring of integers in the 
unramified quadratic extension $\Q_{p^2}$ of $\Q_p$. Fixing an embedding 
of $\Q_{p^2}$ into $\Bbb B$, we have $\Z_{p^2}\hookrightarrow \OB$, and 
we can choose an element $\Pi\in \OB$ with $\Pi^2=p$ 
such that $\Pi a = a^\s\Pi$, for all $a\in \Q_{p^2}$, where 
$\s$ is the generator of the Galois group of $\Q_{p^2}$ over $\Q_p$. 
Then $\OB = \Z_{p^2}[\Pi]$. 

Recall that $W=W(\bar\Bbb F_p)$. Let $\roman{Nilp}$ be the category 
of $W$--schemes $S$ such that $p$ is locally nilpotent in $\Cal O_S$, 
and for $S\in\roman{Nilp}$, let
$\bar S = S\times_W\bar\Bbb F_p$. 

A {\it special formal} (s.f.) $\OB$--module over a $W$ scheme $S$ 
is a p-divisible formal group $X$ over $S$ of dimension $2$ and height $4$
with an action $\iota:\OB\hookrightarrow \End_S(X)$. The Lie algebra $\roman{Lie}(X)$,
which 
is a $\Z_{p^2}\tt\Cal O_S$--module, is required to be free of rank $1$ 
locally on $S$.

Fix a s.f. $\OB$--module $\Bbb X$ over $\Spec(\bar\Bbb F_p)$. Such a 
module is unique up to $\OB$--linear isogeny, and
$$\End^0_{\OB}(\Bbb X) \simeq M_2(\Q_p).\tag4.1$$
We fix such an isomorphism.  
Consider the functor 
$$\Cal D^\bullet: \roman{Nilp} \rightarrow \roman{Sets}\tag4.2$$ 
which associates to each $S\in \roman{Nilp}$ the set of isomorphism 
classes of pairs $(X,\rho)$, where $X$ is a s.f. $\OB$--module over $S$ 
and 
$$\rho: \Bbb X\times_{\Spec(\bar\Bbb F_p)}\bar S \lra X\times_S\bar S\tag4.3$$
is a quasi-isogeny\footnote{This means that, locally on $S$, there is an integer $r$ 
such that $p^r\rho$ is an isogeny.}. The group $\GL_2(\Q_p)$ acts on 
$\Cal D^\bullet$ by 
$$g:(X,\rho) \mapsto (X,\rho\circ g^{-1}).\tag4.4$$ 
There is a decomposition
$$\Cal D^\bullet = \coprod_{i}\Cal D^i,\tag4.5$$
where the isomorphism class of $(X,\rho)$ lies in $\Cal D^i(S)$ 
if $\rho$ has height $i$. The action of $g\in \GL_2(\Q_p)$ 
carries $\Cal D^i$ to $\Cal D^{i-\ord\det(g)}$. 
Drinfeld showed that $\Cal D^\bullet$ is representable 
by a formal scheme, which we also denote by $\Cal D^\bullet$,
and that there is an isomorphism of formal schemes
$$\Cal D^\bullet \isoarrow \hat\O_W^\bullet\tag4.6$$
which is equivariant for the action of $\GL_2(\Q_p)$. 

Similarly, using the notation, $B$, $O_B$, $H$, etc.,  of section 1, the formal scheme 
$\widehat{\Cal M}_p$ represents the functor on $\roman{Nilp}$ 
which associates to $S$ the set of isomorphism classes of triples
$(A,\iota,\bar\eta)$ where $A$ is an abelian scheme of relative dimension 
$2$ over $S$, up to prime to $p$--isogeny, with an 
action $O_B\hookrightarrow \End_S(A)$ satisfying the determinant 
condition (1.7), and $\bar\eta$ is $K^p$--equivalence class of 
$O_B$--equivariant isomorphisms 
$$\eta:\hat V^p(A)\isoarrow B(\A_f^p),\tag4.7$$
where
$$\hat V^p(A) \isoarrow \prod_{\ell\ne p} T_\ell(A)\tt\Q\tag4.8$$ 
is the rational Tate module of $A$. Two isomorphisms
$\eta$ and $\eta'$ are equivalent iff there exists an element $k\in K^p$
such that $\eta'=r(k)\circ\eta$. 

Fix a base point $(A_0,\iota_0,\bar\eta_0)$ in $\widehat{\Cal M}_p(\bar\Bbb F_p)$,
and let
$$\widehat{\Cal M}_p^{\, \sim}:\roman{Nilp}\lra \roman{Sets}\tag4.9$$
be the functor which associates to $S$ the set of isomorphism classes of 
tuples $(A,\iota,\bar\eta,\psi)$, where $(A,\iota,\bar\eta)$ is as before, 
and 
$$\psi:A_0\times_{\Spec(\bar\Bbb F_p)}\bar S \lra A\times_S\bar S\tag4.10$$
is an $O_B$--equivariant p-primary isogeny.  

To relate the functors just defined, we let
$$B' = \End^0(A_0,\iota_0),\qquad\text{and}\qquad H'= (B')^\times,\tag4.11$$ 
and fix 
$\eta_0\in \bar\eta_0$. Since the endomorphisms of $\hat V^p(A_0)$
coming from $B'(\A_f^p) = B'\tt_\Q\hat\Z^p$ commute with $O_B$, the corresponding 
endomorphisms of $B(\A_f^p)$, obtained via $\eta_0$, are given 
by right multiplications by elements of $B(\A_f^p)$.  Thus we obtain 
identifications, as in (1.15), 
$$B'(\A_f^p) \isoarrow B(\A_f^p)^{\roman{op}},\qquad\text{and}\qquad H'(\A_f^p) \isoarrow
H(\A_f^p)^{\roman{op}},\tag4.12$$ 
where the order of multiplication is reversed in $B^{\roman{op}}$ and $H^{\roman{op}}$. 
We also identify $\Bbb B$ with $B_p$ and take $\Bbb X = A_0(p)$, 
the p-divisible group of $A_0$, with the action of $\OB = O_B\tt_\Z\Z_p$
coming from $\iota_0$. Note that we also obtain an identification
$B'_p\simeq \GL_2(\Q_p)$, via (4.1). 
 
Once these identifications have been made, there is a natural isomorphism
$$\widehat{\Cal M}_p^{\, \sim} \isoarrow \Cal D^\bullet\times H(\A_f^p)/K^p\tag4.13$$
defined as follows. To a given $(A,\iota,\bar\eta,\psi)$ over $S$, we 
associate:
$$\align
X &= A(p) = \text{ the p-divisible group of $A$,}\\
\nass
\iota&=\text{ the action of $\OB = O_B\tt_\Z\Z_p$ on $A(p)$,}\tag4.14\\
\nass
\rho&=\rho(\psi) = \text{the quasi-isogeny }\\
\nass
{}&\qquad\qquad \rho(\psi):\Bbb X\times_{\Spec(\bar\Bbb F_p)} \bar S \lra A(p)\times_S\bar S\\
\nass
{}&\qquad\text{determined by $\psi$,}
\endalign
$$
so that $(X,\rho)$ defines an element of $\Cal D^\bullet$.  For $\eta\in \bar\eta$, 
there is
also  a diagram
$$\matrix \hat V^p(A_0)&\overset{\eta_0}\to{\isoarrow}& \hskip -12pt B(\A_f^p)\\
\nass
\nass
\psi_*\downarrow\phantom{\psi_*}&{}&\hskip -12pt\phantom{r(g)}\downarrow r(g)\\
\nass
\hat V^p(A)&\overset{\eta}\to{\isoarrow}&\hskip -12pt B(\A_f^p)
\endmatrix\tag4.15
$$
where $r(g)$ denotes right multiplication\footnote{This is the reason 
for the ``op" in the isomorphism (4.12), since we ultimately 
want a left action of $H'(\A_f^p)$.} by an element $g\in H(\A_f^p)$.
The coset $gK^p$ is then determined by the equivalence class $\bar\eta$, 
and the isomorphism (4.13) sends $(A,\iota,\bar\eta,\psi)$ to 
$((X,\rho), gK^p)$.  

The Drinfeld-Cherednik Theorem says that, by passing to the 
quotient under the action of $H'(\Q)$, we have
$$\matrix 
\widehat{\Cal M}_p^{\, \sim} &\isoarrow& \Cal D^\bullet\times H(\A_f^p)/K^p\\
\nass
\downarrow&{}&\downarrow\\
\nass
\widehat{\Cal M}_p&\isoarrow& H'(\Q)\back\bigg(\Cal D^\bullet\times H(\A_f^p)/K^p\bigg).
\endmatrix\tag4.16
$$
Via the isomorphism $\Cal D^\bullet \overset{\sim}\to{\rightarrow}
\hat\O_W^\bullet$ of (4.6), this yields (1.16). 

We can now describe the formal scheme 
determined by the cycle $\ZZ(t)$, following section 8 of 
\cite{\krinvent}. Let\footnote{We change to $\widehat{\Cal C}$ to avoid confusion 
with the notation $\Zh(t,v)$ used for classes in the arithmetic 
Chow group.}
$\widehat{\Cal C}_p(t)$ be the base change to $W$ of the formal 
completion of $\ZZ(t)$ along its fiber at $p$, and let
$\widehat{\Cal C}_p^{\,\sim}(t)$ be the fiber product:
$$\matrix \widehat{\Cal C}_p^{\,\sim}(t) & \lra &\widehat{\Cal M}_p^{\,\sim}\\
\nass
\downarrow&{}&\downarrow\\
\nass
\widehat{\Cal C}_p(t)&\lra&\widehat{\Cal M}_p.\endmatrix\tag4.17$$
A point of $\widehat{\Cal C}_p(t)$ corresponds to a collection 
$(A,\iota,\bar\eta,x)$, where $x\in V(A,\iota)$ 
is a special endomorphism with
$Q(x)=t$.
In addition, since we are now working with $A$'s up to prime to 
$p$ isogeny, we also require that the endomorphism $\eta_*(x)$ 
of $B(\A_f^p)$, obtained by transfering, via $\eta$, the endomorphism 
of $\hat V^p(A)$ induced by $x$, is given by right
multiplication by an element
$j_f^p(x)\in V(\A_f^p)\cap \hat{O}_B$. This condition does not depend 
on the choice of $\eta$ in the $K^p$--equivalence class $\bar\eta$. 

Next, we would like to determine the image of $\widehat{\Cal C}_p^{\,\sim}(t)$
in $\Cal D^\bullet\times H(\A_f^p)/K^p$ under the isomorphism in the top line of 
(4.16). 
Let
$$V'=\{\ x\in B'\ \mid\ \tr(x)=0\ \} = V(A_0,\iota_0)\tt_\Z\Q,\tag4.18$$
where $V(A_0,\iota_0)$ is the space of special endomorphisms of $(A_0,\iota_0)$. 
We again write $Q$ for the quadratic form on this space. 
For a point of $\widehat{\Cal C}_p^{\,\sim}(t)$
associated to a collection $(A,\iota,\bar\eta,x,\psi)$, let $((X,\rho),gK^p)$ 
be the corresponding point in $\Cal D^\bullet\times H(\A_f^p)/K^p$. 
The special endomorphism $x\in V(A,\iota)$ induces an endomorphism of $A\times_S\bar S$
and, thus, via the isogeny $\psi$, 
there is an associated endomorphism $\psi^*(x)\in V'(\Q)$. This element 
satisfies two compatibility conditions with the other data:
\roster
\item"{(i)}" $\psi^*(x)$ determines an element $j=j(x)\in V'(\Q_p) = \End_{\OB}^0(\Bbb X)$.
By construction, this element has the property that the corresponding element
$$\rho\circ j \circ \rho^{-1} \in \End_{\OB}^0(X\times_S\bar S)\tag4.19$$
is, in fact, the restriction of an element of $\End_{\OB}(X)$, viz. the 
endomorphism induced by original $x$. 
Said another way, 
$j$ defines an endomorphism of the reduction $X\times_S\bar S$ which 
lifts to an endomorphism of $X$.  
\item"{(ii)}" Via the diagram (4.15),  
$$g^{-1}\,\psi^*(x)\,g \in V(\A_f^p)\cap \hat{O}_B^p.\tag4.20$$
\endroster
Here we are slightly abusing notation and, in effect, are 
identifying $\psi^*(x)\in V'(\Q)$ with an element of $V(\A_f^p)$ 
obtained via the identification
$$V'(\A_f^p)\isoarrow V(\A_f^p)\tag4.21$$
coming from (4.12).

Condition (i) motivates the following basic definition, \cite{\krinvent}, 
Definition~2.1.
\proclaim{Definition 4.1} For a special endomorphism $j\in V'(\Q_p)$ of $\Bbb X$, let 
$Z^\bullet(j)$ be the closed formal subscheme of $\Cal D^\bullet$ 
consisting of the points $(X,\rho)$ such that $\rho\circ j\circ \rho^{-1}$ 
lifts to an endomorphism of $X$. 
\endproclaim
We will also write $Z(j)\subset \Cal D$ for the subschemes where the 
height of the quasi-isogeny $\rho$ is $0$. We will give a more detailed 
description of $Z(j)$ in a moment. 

As explained above and in more detail in \cite{\krinvent}, section 8, there is a 
map
$$\widehat{\Cal C}_p^{\,\sim}(t)\hookrightarrow 
V'(\Q)\times\Cal D^\bullet\times H(\A_f^p)/K^p\tag4.22$$
whose image is the set:
$$(\star):=\left\{\ (y,(X,\rho),gK^p)\ \bigg\vert\ \matrix (i)\quad Q(y)=t\\
\nass
 (ii)\quad(X,\rho)\in Z^\bullet(j(y))\\
\nass 
(iii)\quad y\in g\,(V(\A_f^p)\cap \hat{O}_B^p)\,g^{-1}\endmatrix\ \right\}.\tag4.23$$ 
Taking the quotient by the group $H'(\Q)$, we obtain
the following p-adic uniformization of the special cycle:
\proclaim{Proposition 4.2} The construction aboves yields isomorphisms:
$$\matrix
\widehat{\Cal C}_p(t) & \isoarrow & H'(\Q)\back (\star)\\
\nass
\downarrow&{}&\downarrow\\
\nass 
\widehat{\Cal M}_p & \isoarrow & H'(\Q)\back \bigg(\ \Cal D^\bullet\times H(\A_f^p)/K^p\ \bigg)
\endmatrix
$$
of formal schemes over $W$.
\endproclaim
{\bf Remark 4.3.} 
In fact, in this discussion, the compact open subgroup $K^p$ giving the 
level structure away from $p$ can be arbitrary. In the case of interest, where
$K^p = (\hat{O}^p_B)^\times$, the last diagram can be simplified as follows. 
Let
$$L' = V'(\Q) \cap \big(\,B'(\Q_p)\times \hat O_B^p\,\big),\tag4.24$$
so that $L'$ is a $\Z[p^{-1}]$--lattice in $V'(\Q)$. Also, as in section 1, let
$$\Gamma' = H'(\Q)\cap \big(\,H'(\Q_p)\times K^p\,\big),\tag4.25$$
so that (the projection of) $\Gamma'$ is an arithmetic subgroup of $H'(\Q_p)\simeq \GL_2(\Q_p)$. 
Finally, let
$$\Cal D^\bullet_t = 
\coprod_{\matrix\scr y\in L'\\ \scr Q(y)=t\endmatrix} Z^\bullet(j(y)).\tag4.26$$
Then, 
$$\matrix
\widehat{\Cal C}_p(t) & \isoarrow &\big[\,\Gamma'\back \Cal D^\bullet_t\,\big]\\ 
\nass
\downarrow&{}&\downarrow\\
\nass 
\widehat{\Cal M}_p&\isoarrow& \big[\,\Gamma'\back \Cal D^\bullet\,\big].
\endmatrix\tag4.27
$$
Of course, we should now view these quotients as orbifolds, and, in fact, should 
formulate the discussion above in terms of stacks. 

To complete the picture of the vertical components of our cycle $\ZZ(t)$, 
we need a more precise description of the formal schemes $Z^\bullet(j)$ for $j\in V'(\Q_p)$, 
as obtained in the first four sections of \cite{\krinvent}. Note that we are using 
the quadratic form $Q(j) = \det(j)$, whereas, in \cite{\krinvent}, 
the quadratic form $q(j) =j^2=-Q(j)$ was used. 
It is most
convenient  to give this description in the space $\hat\O_W^\bullet=\hat\O_W\times\Z$. 
In fact, we will just consider $Z(j)\subset \hat\O_W$, and we will 
assume that $p\ne2$. The results for $p=2$, which are very similar, are 
described in section 11 of
\cite{\kryII}.  Recall that the 
special fiber of $\hat\O_W$ is a union of projective lines $\Bbb P_{[\L]}$, 
indexed by the vertices $[\L]$ of the building $\Cal B$ of 
$\PGL_2(\Q_p)$. 

The first result, proved in section 2 of \cite{\krinvent}, describes the 
support of $Z(j)$ in terms of the building. 
\proclaim{Proposition 4.4} (i) 
$$\Bbb P_{[\L]}\cap Z(j)\ne\emptyset \iff j(\L)\subset\L.$$
In particular, if $Z(j)\ne\emptyset$, then $\ord_p(Q(j))\ge0$.\hfb
(ii) 
$$j(\L)\subset\L\quad\iff\quad d([\L],\Cal B^j)\ \le\  \frac12\cdot\ord_p(Q(j)).$$
Here $\Cal B^j$ is the fixed point set of $j$ on $\Cal B$, and $d(x,y)$ 
is the distance between the points $x$ and $y\in \Cal B$. 
\endproclaim

Recall that the distance function on the building $\Cal B$ is $\roman{PGL}_2(\Q_p)$--invariant
and gives each edge length $1$. 

In effect, if we write 
$Q(j) = \e\,p^\a$, for $\e\in \Z_p^\times$, then
the support of $Z(j)$ lies in the set of $\Bbb
P_{[\L]}$'s  indexed by vertices of $\Cal B$ in the `tube' $\Cal T(j)$ of radius
$\frac{\a}2$ around the fixed point set $\Cal B^j$ of $j$.  

Next, the following observation of Genestier is essential, \cite{\krinvent}, Theorem~3.1:
\hfb
Write $Q(j) = \e\,p^\a$. Then
$$Z(j) = \cases (\hat\O_W)^j&\text{ if $\a=0$,}\\
\nass
(\hat\O_W)^{1+j}&\text{ if $\a>0$.}
\endcases\tag4.28
$$
where $(\hat\O_W)^x$ denotes the fixed point set of the 
element $x\in \GL_2(\Q_p)$ acting on $\hat\O_W$, via its projection to $\PGL_2(\Q_p)$.
Using this fact, one can obtain local equations for $Z(j)$ 
in terms of the local coordinates on $\hat\O_W$, cf. \cite{\krinvent}, section 3. 
Recall that there are standard coordinate neighborhoods associated to each 
vertex $[\L]$ and each edge $[\L_0,\L_1]$ of $\Cal B$, cf. \cite{\krinvent}, 
section 1. By Proposition~4.4, it suffices to compute in the neighborhoods 
of those 
$\Bbb P_{[\L]}$'s which meet the support of $Z(j)$. For the full result, 
see Propositions~3.2 and~3.3 in
\cite{\krinvent}.  
In particular, embedded components can occur. But since these turn out to be
negligible,  e.g., for intersection theory, cf. section 4 of \cite{\krinvent}, we 
can omit them from $Z(j)$ to obtain the divisor $Z(j)^{\roman{pure}}$
which has the following description. 

\proclaim{Proposition 4.5} Write $Q(j) = \e\,p^\a$, and let
$$\mu_{[\L]}(j) = \max\{\,0, \,\frac{\a}2-d([\L],\Cal B^j)\,\}.$$
(i) If $\a$ is even and $-\e\in \Z_p^{\times, 2}$, then
$$Z(j)^{\roman{pure}} = \sum_{[\L]} \mu_{[\L]}(j)\cdot \Bbb P_{[\L]}.$$
(ii) If $\a$ is even and $-\e\notin \Z_p^{\times, 2}$, then
$$Z(j)^{\roman{pure}} = Z(j)^h+\sum_{[\L]} \mu_{[\L]}(j)\cdot \Bbb P_{[\L]},$$
where the horizontal part $Z(j)^h$ is the disjoint union 
of two divisors projecting isomorphically to $\roman{Spf}(W)$
and meeting the special fiber in `ordinary special' points of $\Bbb P_{[\L(j)]}$, 
where $[\L(j)]$ is the unique vertex $\Cal B^j$ of $\Cal B$ fixed by $j$. \hfb
(iii) If $\a$ is odd, then
$$Z(j)^{\roman{pure}} = Z(j)^h+\sum_{[\L]} \mu_{[\L]}(j)\cdot \Bbb P_{[\L]},$$
where the horizontal divisor $Z(j)^h$ is $\Spf(W')$ where $W'$ is the ring of 
integers in a ramified quadratic extension of $W\tt_{\Z}\Q$ and $Z(j)^h$ meets 
the special fiber in the double point $\roman{pt}_{\Delta(j)}$ 
where $\Delta(j)$ is the edge of $\Cal B$ containing the unique fixed point $\Cal B^j$ 
of $j$. \hfb
Thus $Z(j)^{\roman{pure}}$ is a sum with multiplicities of regular one dimensional 
formal schemes.  
\endproclaim

In case (i), the split case, 
$\Q_p(j) \simeq \Q_p\oplus \Q_p$, the element $j$ lies in a split torus $A$
in $\GL_2(\Q_p)$, and $\Cal B^j$ is the corresponding apartment in $\Cal B$. 
More concretely, if $e_0$ and $e_1$ are eigenvectors of $j$ giving a basis of $\Q_p^2$, then
$\Cal B^j$ is the geodesic arc connecting the vertices $[\Z_p\, e_0 \oplus p^r\Z_p\, e_1]$, for $r\in
\Z$.  The $\Bbb P_{[\L]}$'s for $[\L]\in \Cal B^j$ have multiplicity $\frac{\a}2$ 
in $Z(j)$, and the multiplicity decreases linearly with the distance from $\Cal B^j$. 
The cycle $Z(j)$ is infinite and there are no horizontal components. 

In case (ii), the inert case, $\Q_p(j)$ is an unramified quadratic extension 
of $\Q_p$, the element $j$ lies in the corresponding nonsplit Cartan subgroup 
of $\GL_2(\Q_p)$ and $\Cal B^j = [\L(j)]$ is a single vertex. The corresponding
$\Bbb P_{[\L(j)]}$ occurs with multiplicity $\frac{\a}2$ in $Z(j)$, and the multiplicity
of the vertical components $\Bbb P_{[\L]}$ decreases linearly with 
the distance $d([\L],[\L(j)])$. 

Finally, in case (iii), the ramified case, $\Q_p(j)$ is a ramified 
quadratic extension of $\Q_p$, the element $j$ lies in the corresponding nonsplit Cartan subgroup 
of $\GL_2(\Q_p)$ and $\Cal B^j$ is the midpoint of a unique edge $\Delta(j) = [\L_0,\L_1]$. 
The vertical components $\Bbb P_{[\L_0]}$ and $\Bbb P_{[\L_1]}$ 
occur with multiplicity $\frac{\a-1}2$ in $Z(j)$, and, again, the multiplicity
of the vertical components $\Bbb P_{[\L]}$ decreases linearly with 
the distance $d([\L],\Cal B^j)$, which is now a half-integer. 

This description of the $Z(j)$'s, together with the p-adic unformization 
of Proposition~4.2, gives a fairly complete picture of the vertical components of 
the cycles $\ZZ(t)$ in the fibers $\Cal M_p$, $p\mid D(B)$, of bad reduction. 
Several interesting features are evident. 

For example, the following result gives a
criterion  for the occurrence of such components.
\proclaim{Proposition 4.6} For $p\mid D(B)$, the cycle $\Cal Z(t)$ contains 
components of the fiber $\Cal M_p$ of bad reduction if and only if 
$\ord_p(t)\ge 2$, and no prime $\ell\ne p$ with $\ell \mid D(B)$
is split in $\kay_t$. 
\endproclaim
Note that the condition amounts to (i) $\ord_p(t)\ge2$, and (ii) the field 
$\kay_t$ embeds into $B^{(p)}$

For example, if more than one prime $p\mid D(B)$ 
splits in the quadratic field $\kay_t=\Q(\sqrt{-t})$, then $\ZZ(t)$ 
is empty. If $p\mid D(B)$ splits in $\kay_t$, and all other primes $\ell\mid D(B)$ 
are not split in $\kay_t$, then the generic fiber $\Cal Z(t)_\Q$ is empty, and $\Cal Z(t)$ is 
a vertical cycle in the fiber at $p$, provided $\ord_p(t)\ge2$. 
In general, if $p\mid D(B)$, and $\kay_t$ embeds in $B^{(p)}$, then the vertical component 
in $\Cal M_p$ of the cycle $\Cal Z(p^{2r}t)$ grows as $r$ goes to infinity, while the 
horizontal part does not change.
Indeed, if we change $t$ to $p^{2r}t$, then, by (3.8), 
$\deg_\Q \ZZ(p^{2r}t) = \deg_\Q \ZZ(t)$, while
both the radius of the tube $\Cal T(p^rj)$ and the multiplicity function 
$\mu_{[\L]}(p^rj)$ increase.

Analogues of these results about 
cycles defined by special endomorphisms are obtained in \cite{\krHB} and \cite{\krsiegel}
for Hilbert--Blumenthal varieties and Siegel modular varieties of genus $2$ respectively. 


\subheading{\Sec5. Green functions}

To obtain classes in the arithmetic Chow group $\CH^1_\R(\Cal M)$ from the $\Cal Z(t)$'s,  
it is necessary to equip them with Green functions. These are defined as follows; 
see  
\cite{\annals} for more details. 
For $x\in V(\R)$, with $Q(x)\ne 0$, let 
$$D_x=\{\ z\in D\ \mid\ (x,w(z))=0\ \}.\tag5.1$$
Here $w(z)\in V(\C)$ is any vector with image $z$ in $\Bbb P(V(\C))$.  
The set $D_x$ consists of two points if $Q(x)>0$, and is empty if $Q(x)<0$. 
By (3.5) and (3.6), we have
$$\Cal Z(t)(\C) = \sum_{\matrix\scr x\in O_B\cap V\\ \scr Q(x)=t\\ \scr\mod \Gamma
\endmatrix} \pr(D_x)\tag5.2$$
where $\pr:D \rightarrow \Gamma\back D$. 
For $x\in V(\R)$, with $Q(x)\ne 0$, and $z\in D$, let
$$R(x,z) = |(x,w(z))|^2|(w(z),\overline{w(z)})|^{-1}.\tag5.3$$
This function on $D$ vanishes precisely on $D_x$. Let
$$\beta_1(r) = \int_1^\infty e^{-ru}\,u^{-1}\,du = - \Ei(-r)\tag5.4$$ 
be the exponential integral. Note that 
$$
\beta_1(r) =
\cases -\log(r) -\gamma + O(r),& \text{ as $r\rightarrow 0$,}\\
\nass
O(e^{-r}), &\text{ as $r\rightarrow \infty$.}
\endcases\tag5.5
$$
Thus, the function
$$\xi(x,z) = \beta_1(2\pi R(x,z))\tag5.6$$
has a logarithmic singularity on $D_x$ and decays exponentially 
as $z$ goes to the boundary of $D$. A straightforward calculation, \cite{\annals}, 
section 11, 
shows that $\xi(x,\cdot)$ is a Green function for $D_x$. 
\proclaim{Proposition 5.1} As currents on $D$, 
$$dd^c\xi(x) + \delta_{D_x} = [\ph^0_\infty(x)\,\mu],$$
where, for $z\in D\simeq \C\setminus\R$ with $y = \roman{Im}(z)$, 
$$\mu =\frac1{2\pi}\,\frac{i}2\,\frac{dz\wedge d\bar z}{y^2}$$
is the hyperbolic volume form and 
$$\ph^0_\infty(x,z) = \big[\ 4\pi(\, R(x,z) + 2 Q(x) \, ) -1\ \big]\cdot e^{-2\pi R(x,z)}.$$
\endproclaim
Recall that $dd^c = -\frac{1}{2\pi i} \partial\db$. 

{\bf Remark.} For fixed $z\in D$, 
$$(x,x)_z = (x,x) + 2\,R(x,z)\tag5.7$$
is the majorant attached to $z$, \cite{\siegel}. Thus, the function
$$\ph_\infty(x,z)\cdot \mu = \ph_\infty^0(x,z)\cdot e^{-2\pi Q(x)}\cdot \mu\tag5.8$$
is (a very special case of) the Schwartz function valued in smooth 
$(1,1)$--forms on $D$ defined in 
\cite{\kmI}.  In fact, the function $\xi(x,z)$ was first obtained 
by solving the Green equation of Proposition~5.1  with this right hand side.

Because of the rapid decay of $\xi(x,\cdot)$, we can average over lattice points.
\proclaim{Corollary 5.2} For $v\in \R^\times_{>0}$, let
$$\Xi(t,v)(z) := \sum_{\scr x \in O_B\cap V\atop \scr Q(x)=t} \xi(v^\frac12 x,z).$$
(i) For $t>0$, 
$\Xi(t,v)$ defines a Green function for $\Cal Z(t)$. \hfb
(ii) For $t<0$, $\Xi(t,v)$ defines a smooth function on $M(\C)$. 
\endproclaim
Note, for example, that for $t>0$, 
\vskip -10pt
$$\matrix \Cal Z(t)(\C)=\emptyset\\ \nass \text{ and }\Xi(t,v)=0
\endmatrix\ \iff \ \{ x \in O_B\cap V,\   Q(x)=t\} =
\emptyset\ \iff\ 
\matrix\text{$\kay_t$  does not}\\\text{ embed in $B$.}\endmatrix \tag5.9
$$

An explicit construction of Green functions for divisors in general locally 
symmetric varieties is given by Oda and Tsuzuki, \cite{\odatsuzuki}, by a different 
method.

\subheading{\Sec6. The arithmetic theta series}

At this point, we can define a family of classes in $\CH^1(\Cal M)$. These 
can be viewed as an analogue for the arithmetic surface $\Cal M$ 
of the Hirzebruch-Zagier classes $T_N$ in the middle cohomology of 
a Hilbert modular surface, \cite{\hirzebruchzagier}. 

\proclaim{Definition 6.1} For $t\in \Z$, with $t\ne0$, and for a parameter 
$v\in \R^\times_+$, define classes in $\CH^1_\R(\Cal M)$ by
$$\Zh(t,v) = \cases (\,\Cal Z(t), \Xi(t,v)\,)&\text{if $t>0$,}\\
\nass
(\,0,\Xi(t,v)\,)&\text{if $t<0$.}
\endcases
$$
For $t=0$, define
$$\Zh(0,v) = - \ohat -(0,\log(v)) + (0,\bold c),$$
where $\bold c$ is the constant defined by (2.15). 
\endproclaim

We next construct a generating series for these classes; 
again, this can be viewed as an arithmetic analogue of the Hirzebruch--Zagier 
generating series for the $T_N$'s. For $\tau=u+iv\in \H$, the 
upper half plane, let $q=e(\tau) = e^{2\pi i\tau}$. 

\proclaim{Definition 6.2} The arithmetic theta series is the generating series
$$\phih(\tau) = \sum_{t\,\in\, \Z} \Zh(t,v)\,q^t\ \in \CH^1(\Cal M)[[q]].$$
\endproclaim
Note that, since the imaginary part $v$ of $\tau$ appears as a parameter in the 
coefficient $\Zh(t,v)$, this series is not a holomorphic function of $\tau$. 

The arithmetic theta function $\phih(\tau)$ is closely connected with the 
generating series for quadratic divisors considered by Borcherds, \cite{\borchduke}, and 
the one for Heegner points introduced by Zagier,
\cite{\zagier}.  
The following result, \cite{\ariththeta}, justifies the terminology. 
Its proof, which will be sketched in section 7, depends on 
the results of \cite{\kryII} and of Borcherds, \cite{\borchduke}.

\proclaim{Theorem 6.3} The arithmetic theta series $\phih(\tau)$ is a 
(nonholomorphic) modular form of weight $\frac32$  
valued in $\CH^1(\Cal M)$. 
\endproclaim

As explained in Proposition~7.1 below, the arithmetic Chow group $\CH^1(\Cal M)$
can be written as direct sum 
$$\CH^1(\Cal M) = \CH^1(\Cal M,\mu)_\C\oplus C^\infty_0(\Cal M(\C))$$
where $\CH^1(\Cal M,\mu)_\C$ a finite dimensional complex vector space, the
Arakelov Chow group of $\Cal M$ for the hyperbolic metric $\mu$, (7.7), and 
$C^\infty_0(\Cal M(\C))$ is the space of smooth functions on $\Cal M(\C)$ 
with integral $0$ with respect to $\mu$. Theorem~6.3 then means that there 
is a smooth function of $\phi_{\roman{Ar}}$ of $\tau$ valued in $\CH^1(\Cal M,\mu)_\C$
and a smooth function $\phi(\tau,z)$ on $\H\times \Cal M(\C)$, with 
$$\int_{\Cal M(\C)}\phi(\tau,z)\,d\mu(z)=0,$$
and such that the sum $\phi(\tau)=\phi_{\roman{Ar}}(\tau)+\phi(\tau,z)$ satisfies the 
usual transformation law for a modular form of weight $\frac32$, 
for $\Gamma_0(4 D(B))$, and such that the $q$--expansion 
of $\phi(\tau)$ is the formal generating series $\thh(\tau)$ of 
Definition~6.2. By abuse of notation, we write $\thh(\tau)$ both for $\phi(\tau)$
and for its $q$--expansion. 

\subheading{\Sec7. Modularity of the arithmetic theta series}

Although the arithmetic theta series $\phih(\tau)$ can be viewed as a kind 
of generating series for lattice vectors in the spaces $V(A,\iota)$ 
of special endomorphisms, there is no evident analogue of the Poisson 
summation formula, which is the key ingredient of the proof of 
modularity of classical theta series. Instead, the modularity of $\phih(\tau)$ 
is proved by computing its height pairing with generators of the group $\CH^1(\Cal M)$ 
and identifying the resulting functions of $\tau$ with known modular forms.  

First we recall the structure of the arithmetic Chow group $\CH^1_\R(\Cal M)$, 
\cite{\gsihes}, \cite{\soulebook}, \cite{\bost}. 
There is a map
$$a: C^\infty(M(\C)) \lra \CH^1_\R(\Cal M), \qquad \phi\longmapsto (0,\phi)\tag7.1$$
from the space of smooth functions on the curve $M(\C)$. Let 
$\triv = a(1)$ be the image of the constant function. Let 
$\Ver \subset \CH^1_\R(\Cal M)$ be the subspace generated by classes of the 
form $(Y_p,0)$, where $Y_p$ is a component of a fiber $\Cal M_p$. The 
relation $\divh(p)=(\Cal M_p,-\log(p)^2)\equiv 0$ implies that 
$(\Cal M_p,0) \equiv 2\log(p)\cdot\triv$, 
so that $\triv\in \Ver$. The spaces $a(C^\infty(M(\C))$ and $\Ver$ span the 
kernel of the restriction map
$$\roman{res}_\Q: \CH^1_\R(\Cal M) \lra \roman{CH}^1(\Cal M_\Q)_\R\tag7.2$$
to the generic fiber. Here $\roman{CH}^1(\Cal M_\Q)_\R=\roman{CH}^1(\Cal M_\Q)\tt_\Z\R$.

Let 
$$\MW_\R = \MW(M)_\R = \roman{Jac}(M)(\Q)\tt_\Z\R\tag7.3$$
be the Mordell--Weil space of the Shimura curve $M=\Cal M_\Q$, and recall 
that this space is the kernel of the degree map
$$\MW_\R \lra \roman{CH}^1(M)_\R \overset{\deg}\to{\lra}\  \R.\tag7.4$$ 
Note that the degree of the restriction of the class $\ohat$ to the 
generic fiber is positive, 
since it is given by the integral of the hyperbolic volume form over $[\Gamma\back D]$. 

Finally, let $C^\infty_0 = C^\infty(M(\C))_0$ be the 
space of smooth functions which are orthogonal to the constants with 
respect to the hyperbolic volume form. 

The following is a standard result in the Arakelov theory of 
arithmetic surfaces, \cite{\hriljac}, \cite{\faltings}, \cite{\bost}. 
\proclaim{Proposition 7.1} Let
$$\MWt := \bigg(\ \R\,\ohat \oplus \Ver \oplus a(C^\infty_0)\ \bigg)^\perp$$
be the orthogonal complement of $\R\,\ohat \oplus \Ver \oplus a(C^\infty_0)$ 
with respect to the height pairing. 
Then 
$$\CH^1_\R(\Cal M) = \MWt \oplus \bigg(\ \R\,\ohat\oplus \Ver\ \bigg)\oplus a(C^\infty_0),$$ 
where the three summands are orthogonal with respect to the 
height pairing. Moreover, the restriction map $\roman{res}_\Q$ 
induces an isometry
$$\roman{res}_\Q:\ \MWt \isoarrow \MW$$
with respect to the Gillet--Soul\'e--Arakelov height pairing on $\MWt$ 
and the negative of the Neron--Tate height pairing on $\MW$. 
\endproclaim

In addition, there are some useful formulas for the height pairings 
of certain classes. For example, for any $\Zh\in \CH^1\!(\Cal M)$, 
$$\gs{\Zh}{\triv} = \frac12\,\deg_\Q(\Zh).\tag7.5$$
In particular, $\triv$ is in the radical of the restiction of 
the height pairing to $\Ver$, and $\gs{\ohat}{\triv} = \frac12\,\deg_\Q(\ohat)>0$. 
Also, for $\phi_1$ and $\phi_2\in C^\infty_0$, 
$$\gs{a(\phi_1)}{a(\phi_2)} = \frac12\,\int_{M(\C)} dd^c\phi_1\cdot \phi_2.\tag7.6$$
Note that the subgroup 
$$\CH^1\!(\Cal M,\mu)_\R = \MWt \oplus  \R\ohat\oplus \Ver\tag7.7$$
is the Arakelov Chow group (with real coefficients) for the hyperbolic 
metric $\mu$. 

Returning to the arithmetic theta series $\phih(\tau)$, we consider its pairing
with various classes in $\CH^1_\R(\Cal M)$. 
To describe these, we first introduce an Eisenstein series of 
weight $\frac32$ associated to the quaternion algebra $B$, \cite{\kryII}. 
Let $\Gamma'=\SL_2(\Z)$ 
and let $\Gamma'_\infty$ be the stabilizer of the cusp at infinity. 
For $s\in \C$, let 
$$\Cal E(\tau,s;B) = v^{\frac12(s-\frac12)}\sum_{\gamma\in \Gamma'_\infty\back \Gamma'}
(c\tau+d)^{-\frac32}\,|c\tau+d|^{-(s-\frac12)}\,\Phi^B(\gamma,s),\tag7.8$$
where $\gamma = \pmatrix a&b\\c&d\endpmatrix$ and $\Phi^B(\gamma,s)$ 
is a function of $\gamma$ and $s$ depending on $B$. This series 
converges absolutely for $\Re(s)>1$ and has an analytic continuation 
to the whole $s$--plane. It satisfies the functional equation, \cite{\kryII}, section 16, 
$$\Cal E(\tau,s;B) = \Cal E(\tau,-s;B),\tag7.9$$ 
normalized as in Langlands theory. 

The following result is proved in \cite{\kryII}. 
\proclaim{Theorem 7.2} \hfb
(i)
$$\align
2\cdot\gs{\phih(\tau)}{\triv} &= -\vol(M(\C))+\sum_{t>0} \deg(\Cal Z(t)_\Q)\,q^t \\
\nass
{}&= \Cal E(\tau,\frac12;B).\\
\noalign{(ii)}
\gs{\phih(\tau)}{\ohat} &= \sum_t \gs{\Zh(t,v)}{\ohat}\,q^t \\
\nass
{}&= \Cal E'(\tau,\frac12;B).
\endalign
$$
\endproclaim

This result is proved by a direct computation of the Fourier 
coefficients on the two sides of (i) and (ii). For (i), 
the identity amounts to the formula
$$
\deg(\Cal Z(t)_\Q)\,q^t = 2\,\delta(d;D(B))\,H_0(t;D(B))\,q^t = \Cal
E'_t(\tau,\frac12;D(B)),\tag7.10
$$
with the notation as in (3.8), (3.9) and (3.10) above. Recall that
we write $4t = n^2 d$ where $-d$ is the discriminant of the field $\kay_t$.

For example, if $D(B)=1$, i.e., in the case of $B=M_2(\Q)$, 
$$H_0(t;1) = \sum_{c|n} \frac{h(c^2d)}{w(c^2d)}.\tag7.11$$
Thus, $2\,H_0(t;1) = H(4t)$, where $H(t)$ is the `class number' which appears in 
the Fourier expansion of Zagier's nonholomorphic Eisenstein series of weight $\frac32$, 
\cite{\cohen},
\cite{\zagier}:
$$\Cal F(\tau)= -\frac1{12} + \sum_{t>0} H(t)\,q^t + \sum_{m\in \Z} 
\frac1{16\pi}\, v^{-\frac12}
\int_1^\infty e^{-4\pi m^2 v r}\, r^{-\frac32}\,dr\, q^{-m^2}.\tag7.12$$
In fact, when $D(B) =1$, the value of our Eisenstein series is
$$\Cal E(\tau,\frac12;1)= -\frac1{12} + \sum_{t>0} 2\,H_0(t;1)\,q^t+
\sum_{m\in\Z}\frac{1}{8\pi}\,v^{-\frac12}\,\int_1^\infty e^{-4\pi m^2 vr}\,
r^{-\frac32}\,dr \cdot q^{-m^2}.\tag7.13$$
As in (i) of Theorem~7.2, both of these series have an interpretation in 
term of degrees of $0$--cycles of CM points  on the modular curve, \cite{\yangMSRI}. 
Their relation to a regularized integral of a theta series is proved by J.~Funke 
in \cite{\funkethesis}, \cite{\funkecompo}. 

The computations involved in (ii) are significantly more difficult. 
For example, if $t>0$ and $\Cal Z(t)_\Q(\C)$ is nonempty, then, 
\cite{\kryII}, Theorem~8.8, the $t$-th
Fourier  coefficient of $\Cal E'(\tau,\frac12;B)$ is
$$\align
&\Cal E'_t(\tau,\frac12;B)\tag7.14\\
\nass
{}&= 2\,\delta(d;D)\,H_0(t;D)\cdot q^t\cdot\,\bigg[ \,
\frac12\, \log(d)
 + \frac{L'(1,\chi_d)}{L(1,\chi_d)} -\frac12\log(\pi) -\frac12\gamma\\
\nass
\nass
{}&\hbox to .7in{} 
+\frac12 J(4\pi tv) 
+\sum_{p \atop p\nmid D} \bigg(\ \log|n|_p - \frac{b_p'(n,0;D)}{b_p(n,0;D)}\ \bigg)
+\sum_{p \atop p\mid D} K_p\,\log(p)
\ \bigg].\\ 
\endalign
$$
Here, we write $D$ for $D(B)$, $4t=n^2d$, as before, and 
$$K_p = \cases -k + \frac{(p+1)(p^k-1)}{2(p-1)}&\text{ if $\chi_d(p)=-1$, and}\\
\nass
-1-k + \frac{p^{k+1}-1}{p-1} &\text{ if $\chi_d(p)=0$,}
\endcases\tag7.15
$$
with $k=k_p=\ord_p(n)$. Also,  
$$J(x) = \int_0^\infty e^{- x r}\big[\, (1+r)^{\frac12} - 1\,\big]\,r^{-1}\,dr.\tag7.16$$
Finally, for a prime $p\nmid D$, 
$$
\frac{1}{ \log p}\cdot\frac{b_p'(n, 0; D)}{b_p(n, 0; D)}
=\frac{\chi_d(p) -\chi_d(p) \,(2k+1) p^{k} +(2k+2) p^{k+1})}
       {1-\chi_d(p)+\chi_d(p)\, p^{k} -p^{k+1}} 
  -\frac{2p}{1-p}
\tag7.17
$$
The connection of this rather complicated quantity with arithmetic geometry is 
not at all evident.  
Nonetheless, the identity of part (ii) of Theorem~7.2 asserts that
$$\Cal E'_t(\tau,\frac12;B) = \gs{\Zh(t,v)}{\ohat}\,q^t.\tag7.18$$
Recall that the points of the generic fiber $\Cal Z(t)_\Q$ correspond to 
abelian surfaces $A$ with an action of $O_B\tt_\Z\Z[\sqrt{-t}]$, an order in 
$M_2(\kay_t)$. The contribution to $\gs{\Zh(t,v)}{\ohat}$ of the associated 
horizontal component of $\Cal Z(t)$ is the Faltings height of $A$. 
Due to the action of $O_B\tt_\Z\Z[\sqrt{-t}\,]$, $A$ is isogenous to a 
product $E\times E$ where $E$ is an elliptic curve with CM by the order $O_d$, and so the 
Faltings heights are related by 
$$h_{\roman{Fal}}(A) = 2 \, h_{\roman{Fal}}(E) + \text{an isogeny correction.}\tag7.19$$
Up to some combinatorics involving 
counting of the points $\Cal Z(t)_\Q$, the terms on the first line of the right side of
(7.14) come from $h_{\roman{Fal}}(E)$, while sum on $p\nmid D(B)$ on the second line 
arises from the isogeny correction. The $J(4\pi t v)$ term comes from the 
contribution of the Green function $\Xi(t,v)$ with parameter $v$. Finally, 
the cycle $\Cal Z(t)$ can have vertical components, and the sum on $p\mid D(B)$ 
comes from their pairing with $\ohat$. 

The analogue of (ii) also holds for modular curves, i.e., when $D(B)=1$. In this case, 
there are some additional terms, as in the nonholomorphic parts of (7.12) and (7.13). These 
terms are contributions from a class in $\CH^1_\R(\Cal M)$ supported in the cusp, 
\cite{\yangMSRI}, \cite{\kryIV}. 

To continue the proof of modularity of $\phih(\tau)$, we next consider the 
height pairing of $\phih(\tau)$ with classes of the form $(Y_p,0)$ in $\Ver$, and with
classes of the form $(0,\phi)$, for $\phi\in C^\infty_0$, which 
might be thought of as `vertical at
infinity'. 

First we consider $\gs{\phih(\tau)}{(Y_p,0)}$. 
For $p\ne2$, the intersection number of a component $Y_p$ indexed by $[\L]$ 
with a cycle $\Cal Z(t)$ is calculated explicitly in \cite{\krinvent}. 
For $p=2$, the computation is very similar.
Using this result, we obtain, \cite{\ariththeta}:

\proclaim{Theorem 7.3} Assume that $p\mid D(B)$ and that $p\ne2$. Then, for a 
component $Y_p$ of $\Cal M_p$ associated to a homothety class of lattices $[\Lambda]$, 
there is a Schwartz function $\ph_{[\L]}\in S(V^{(p)}(\A_f))$ and associated 
theta function of weight $\frac32$ 
$$\theta(\tau,\ph_{[\L]}) = \sum_{x\in V^{(p)}(\Q)} \ph_{[\L]}(x)\,q^{Q(x)},$$
such that
$$\gs{\phih(\tau)}{(Y_p,0)} = \theta(\tau,\ph_{[\L]}).$$
\endproclaim

Next consider the height pairings with classes of the form $a(\phi) = (0,\phi)$ 
for $\phi\in C^\infty_0$. Here we note that, for a class $(Z,g_Z)\in \CH^1_R(\Cal M)$, 
$$\gs{(Z,g_Z)}{a(\phi)} = \frac12\,\int_{\Cal M(\C)} \o_Z\,\phi,\tag7.20$$
where $\o_Z = dd^cg_Z + \delta_Z$ is the smooth $(1,1)$--form on the right side of the 
Green equation. The map $(Z,g_Z)\mapsto\o_Z$ defines a map, \cite{\gsihes}, \cite{\soulebook}
$$\o:\CH^1_\R(\Cal M) \lra A^{(1,1)}(\Cal M(\C)).\tag7.21$$
By the basic construction of the Green function
$\Xi(t,v)$ and Proposition~5.1,  we have
$$\align
\o(\,\phih(\tau)\,) & = \sum_{x\in O_B\cap V}
\ph_\infty^0(xv^{\frac12},z)\,q^{Q(x)}\cdot \mu\tag7.22\\
\nass
{}:&=\theta(\tau,\ph^0_\infty)(z)\cdot \mu,
\endalign
$$
where $\theta(\tau,\ph^0_\infty)(z)$ is the 
theta series of weight $\frac32$ for the rational quadratic
space 
$V$ of signature $(1,2)$. 
Thus,
$$\gs{\phih(\tau)}{a(\phi)} = \frac12\,\int_{\Cal M(\C)} 
\theta(\tau,\ph^0_\infty)(z)\,\phi(z)\,d\mu(z),\tag7.23$$
is the classical theta lift of $\phi$. To describe this more 
precisely, recall that
$$dd^c\phi = \frac12\,\Delta\{\phi\}\cdot \mu\tag7.24$$
where $\Delta$ is the hyperbolic Laplacian, and consider functions 
$\phi_\l\in C^\infty_0(\Cal M(\C))$
satisfying $\Delta\phi_\l +\l\phi_\l =0$ for $\l>0$. 
Let $\phi_{\l_i}$, for $0<\l_1\le\l_2\le \dots$  be a basis of such eigenfunctions, 
orthonormal with respect to $\mu$. These are just the Maass forms of 
weight $0$ for the cocompact Fuchsian group $\Gamma = O_B^\times$. 
Recalling (7.6), 
we see that the classes $a(\phi_\l)\in \CH^1_\R(\Cal M)$ are orthogonal with respect to the 
height pairing and span the subspace $a(C^\infty_0)$. By (7.23), we have the 
following, \cite{\ariththeta}. 
\proclaim{Theorem 7.4} For a Maass form $\phi_\l$ of weight $0$ for $\Gamma$, let 
$$\theta(\tau;\phi_\l) :=\int_{\Cal M(\C)} \theta(\tau,\ph^0_\infty)(z)\,\phi_\l(z)\,d\mu(z)$$
be its classical theta lift, a Maass form of weight $\frac32$ and 
level $4D(B)$. Then
$$\gs{\phih(\tau)}{a(\phi_\l)} = \frac12\,\theta(\tau;\phi_\l).$$
\endproclaim

Finally, we must consider the component $\phih_{\MW}(\tau)$ of $\phih(\tau)$ 
in the space $\MWt$. 
Recall the isomorphism
$$\res_\Q:\ \MWt \isoarrow 
\MW\subset\roman{CH}^1(M).\tag7.25$$  
Write
$$\theta_B(\tau):=\res_\Q(\phih(\tau))  = 
-\o+\sum_{t>0} Z(t)\,q^t\quad \in \roman{CH}^1(M),\tag7.26$$
where $Z(t)$ (resp.\ $\o$) is the class of the $0$--cycle $\Cal Z(t)_\Q$ (resp.\ the 
line bundle $\o = \res_\Q(\ohat)$) in $\roman{CH}^1(M)$.
This series is essentially a special case of the generating 
function for divisors considered by 
Borcherds, \cite{\borchduke}, \cite{\borchdukeII}. Assuming that certain spaces of 
vector valued modular forms have bases with rational Fourier coefficients, 
Borcherds proved that his generating series 
are modular forms, and McGraw, \cite{\mcgraw}, verified Borcherds 
assumption. The main point in Borcherds' proof is the existence 
of enough relations among the divisors in question, and such relations 
can be explicitly given via Borcherds construction of meromorphic 
modular forms with product expansions, \cite{\borchinventI}, 
\cite{\borchinventII}. Thus, the series $\theta_B(\tau)$ is a
modular form of weight $\frac32$ valued in $\roman{CH}^1(M)$. 

By (i) of Theorem~7.2, we have
$$\deg(\theta_B(\tau)) =  \Cal E(\tau,\frac12;B),\tag7.27$$
and so, the function
$$\theta_{\MW}(\tau):= \theta_B(\tau) - 
\Cal E(\tau,\frac12;B)\cdot \frac{\o}{\deg(\o)}\quad\in \MW\tag7.28$$
is also modular of weight $\frac32$. Thus we obtain, \cite{\ariththeta}, 
\proclaim{Theorem 7.5} The image of $\phih_{\MW}(\tau)$ 
under the isomorphism (7.25), is given by
$$\align
\res_\Q(\ \phih_{\MW}(\tau)\ ) &= 
\res_\Q\big(\ \phih(\tau) - \Cal E(\tau,\frac12;B)\deg_\Q(\ohat)^{-1}\cdot\ohat\ \big)\\
\nass
{}&= \theta_{\MW}(\tau)\in \MW.
\endalign
$$
Thus, $\phih_{\MW}(\tau)$ is a modular form of weight $\frac32$. 
\endproclaim

This completes the proof of the modularity of the arithmetic theta function 
$\phih(\tau)$. 

\subheading{\Sec8. The arithmetic theta lift}

The arithmetic theta function $\thh(\tau)$ can be 
used to define an arithmetic theta lift
$$\thh:\ S_{\frac32} \lra \CH^1(\Cal M),\qquad f\mapsto \thh(f),$$ 
where $S_{\frac32}$ is the space of cusp forms of weight $\frac32$
for $\Gamma'=\Gamma_0(4D(B))$, as follows. 
Given a cusp form $f$ of weight $\frac32$ for $\Gamma'=\Gamma_0(4D(B))$, let
$$\thh(f) = \gs{\,f}{\thh\,}_{\text{\rm Pet}} 
=\int_{\Gamma'\back \H} f(\tau)\,\overline{\thh(\tau)}\,
v^{\frac32}\,\frac{du\,dv}{v^2}\quad\in \CH^1(\Cal M).\tag8.1$$ 
This construction is analogous to the construction of Niwa, \cite{\niwa}, 
of the classical Shimura lift, \cite{\shimurahalf}, from $S_{\frac32}$ to 
modular forms of weight $2$ for $O_B^\times$. 
Recall that we extend the height pairing $\gs{\ }{\ }$ on $\CH^1_\R(\Cal M)$ 
to an Hermitian pairing on $\CH^1(\Cal M)$, conjugate linear in the 
second argument. By adjointness, and Theorem~7.2,  
$$\langle \thh(f),\triv\rangle = \langle f, \langle
\thh,\triv\rangle\,\rangle_{\text{\rm Pet}} = \frac12\,\langle f,
\Cal E(\frac12;B)\rangle_{\text{\rm Pet}}= 0,\tag8.2$$
and
$$\langle \thh(f),\hat\o\rangle = \langle f, \langle
\thh,\ohat\rangle\,\rangle_{\text{\rm Pet}} = \langle f,
\Cal E'(\frac12;B)\rangle_{\text{\rm Pet}}= 0.\tag8.3$$
If $f$ is holomorphic, then, by Theorem~7.4, for any $\phi\in C^\infty_0(\Cal M(\C))$, 
$$\langle \thh(f),a(\phi)\rangle =\gs{f}{\gs{\thh}{a(\phi)}}
= \gs{f}{\theta(\phi)}_{\roman{Pet}} = 0,\tag8.4$$
since $\theta(\tau;\phi)$ is a combination of Maass forms of weight $\frac32$. 
Thus, for $f$ holomorphic
$$\thh(f)\in \MWt \oplus  \Ver^0\ \subset \CH^1(\Cal M)\tag8.5$$
where 
$$\Ver^0 = \Ver\cap \ker\gs{\cdot }{\ohat}.\tag8.6$$
It remains to describe the components of $\hat\theta(f)$ in the 
spaces $\Ver$ and $\MWt$. 
This is best expressed in terms of automorphic representations.

\subheading{\Sec9. Theta dichotomy: Waldspurger's theory}

We begin with a brief review of Waldspurger's theory of the 
correspondence between cuspidal automorphic representations of the metaplectic cover 
of $\SL_2$ and automorphic representations of $\roman{PGL}_2$ and its 
inner forms.  
For a more detailed survey, the reader can consult
\cite{\pswald}, \cite{\waldsurvey}, \cite{\howeps}, as well as the original 
papers \cite{\waldshimura}, \cite{\waldfourier}, and especially, \cite{\waldquaternion}. 

We fix the additive character $\psi$ of $\A/\Q$ which 
has trivial conductor, i.e., is trivial on $\widehat\Z = \prod_{p<\infty}\Z_p$
and has archimedean component $\psi_\infty(x) = e(x) = e^{2\pi ix}$. 
Since $\psi$ is fixed, we suppress it from the notation. 

Let $G=\SL_2$ and let $G'_\A$ be the $2$--fold cover of $G(\A)$ 
which splits over $G(\Q)$. Let
$$\align
\Cal A_{0}(G')&=\text{the space of genuine cusp forms for $G'_\A$.}\tag9.1\\
\nass
\Cal A_{00}(G')&=\text{the space of genuine cusp forms for $G'_\A$.}\\
&\quad\quad\text{orthogonal to all $O(1)$ theta series}
\endalign
$$
For an irreducible cuspidal automorphic representation $\s\simeq \tt_{p\le\infty}\s_p$ 
in $\Cal A_{00}(G')$, Waldspurger constructs an irreducible cuspidal automorphic representation 
$$\pi = \pi(\s) = \Wald(\s) =\Wald(\s,\psi) \tag9.2$$
of $\roman{PGL}_2(\A)$, which 
serves as a kind of reference point for the description of the global theta lifts
for various ternary quadratic spaces. For each $p\le \infty$, there is a 
corresponding local construction $\s_p\mapsto \Wald(\s_p,\psi_p)$, and the local and 
global constructions are compatible, i.e., 
$$\Wald(\s,\psi)\ \simeq\ \tt_{p\le \infty}\Wald(\s_p,\psi_p).$$
Since we have fixed the additive character $\psi=\tt_p\psi_p$, we will 
often omit it from the notation. 

For a quaternion algebra $B$ over $\Q$, let
$$V^B = \{\ x\in B\mid \tr(x)=0\ \},\tag9.3$$
with quadratic form\footnote{Note that Waldspurger 
uses the opposite sign.} $Q(x) = -x^2=\nu(x)$,
and let
$$H^B = B^\times \simeq \roman{GSpin}(V).\tag9.4$$ 
For a Schwartz function $\ph\in S(V^B(\A))$, $g'\in G'_\A$ and $h\in H^B(\A)$, 
define the theta kernel by
$$\theta(g',h;\ph) = \sum_{x\in V^B(\Q)}\o(g')\ph(h^{-1}x).\tag9.5$$
Here $\o=\o_\psi$ is the Weil representation of $G'_\A$ on $S(V(\A))$.
For $f\in\s$ and $\ph\in S(V^B(\A))$, the classical theta 
lift is
$$\theta(f;\ph) =\gs{f}{\theta(\ph)}_{\roman{Pet}} =
\int_{G'_\Q\back G'_\A}f(g')\,\overline{\theta(g',h;\ph)}\,dg'\quad \in \Cal A_0(H^B).\tag9.6$$
The global
theta lift of $\s$ to $H^B$ is the space
$$\theta(\s;V^B)\ \subset\ \Cal A_0(H^B),\tag9.7$$ 
spanned by the $\theta(f;\ph)$'s for $f\in \s$ and $\ph\in S(V(\A))$. 
Here $\Cal A_0(H^B)$ is the space of cusp forms on $H^B(\Q)\back H^B(\A)$. 
Then $\theta(\s;V^B)$ is either zero or is an irreducible cuspidal representation 
of $H^B(\A)$, \cite{\waldquaternion}, Prop. 20, p. 290.

For any irreducible admissible genuine 
representation $\s_p$ of $G'_p$, the metaplectic cover of 
$G(\Q_p)$, there are analogous local theta lifts
$\theta(\s_p,V_p^B)$.  Each of them is either zero or is an 
irreducible admissible representation of $H^B_p$. The local and global 
theta lifts are compatible in the sense that
$$\theta(\s;V^B) \simeq\cases \tt_p \theta(\s_p;V_p^B)&\text{or}\\
\nass
0.&{}
\endcases\tag9.8
$$
In particular, the global theta lift is zero if any local $\theta(\s_p;V^B_p)=0$, 
but the global theta lift can also vanish even when there is no such local 
obstruction, i.e., even if $\tt_{p\le \infty}\theta(\s_p;V_p^B)\ne0$. 

For each $p$, let $B_p^{\pm}$ be the quaternion algebra over $\Q_p$ 
with invariant 
$$\roman{inv}_p(B^\pm_p) = \pm1.\tag9.9$$ 
Then the two ternary quadratic 
spaces 
$$V_p^{\pm} = \{\ x\in B^\pm_p \mid \tr(x)=0\ \}\tag9.10$$
have the same discriminant and opposite Hasse invariants. 
A key local fact established by Waldspurger is:
\proclaim{Theorem 9.1}{\bf (Local theta dichotomy)} For an irreducible admissible 
genuine representaton $\s_p$ of $G'_p$, precisely 
one of the spaces
$\theta(\s_p;V_p^+)$ and 
$\theta(\s_p;V_p^-)$
is nonzero.
\endproclaim

{\bf Definition 9.2.} Let $\e_p(\s_p)=\pm1$ be the unique sign 
such that with 
$$\theta(\s_p,V_p^{\e_p(\s_p)})\ne0.$$

{\bf Examples:} (i) If $p$ is a finite prime and $\s_p$ 
is an unramified principal series representation, 
then $\e_p(\s_p)=+1$ and $\theta(\s_p;V_p^+)$ 
is a principal series representation of $\GL_2(\Q_p)$. It is unramified if 
$\s_p$ is unramified. \hfb
(ii) For a finite prime $p$ and a character 
$\mu_p$ of $\Q_p^\times$, with $\mu_p^2=|\ |$, there is a 
special representation $\s_p(\mu_p)$ of $G'_p$. For $\s_p=\s_p(\mu_p)$, with 
$\mu_p = |\ |^{\frac12}$, 
$$\align
\theta(\s_p, V_p^-) &=\triv \ne0,\qquad B_p^- = \Bbb B_p \\
\nass
\theta(\s_p, V_p^+) &=0, \qquad\quad\quad B_p^+ = M_2(\Q_p)\\
\nass
\Wald(\s_p, \psi_p)&=\text{ unramified special $\s(|\ |^{\frac12},|\ |^{-\frac12})$ of $GL_2(\Q_p)$.}
\endalign
$$
If $\mu_p\ne|\ |^{\frac12}$, then 
$$\theta(\s_p,V^+) = \Wald(\s_p,\psi_p) = \s(\mu,\mu^{-1})$$
is a special representation of $\GL_2(\Q_p)$ and $\theta(\s_p,V_p^-)=0$. \hfb
(iii) If $\s_\infty = \text{HDS}_{\frac32}$, the holomorphic discrete series 
representation of $G'_\R$ of weight $\frac32$, then
$$\align
\theta(\s_\infty, V_\infty^-) &=\triv \ne0,\qquad B_\infty^- = \Bbb H\\
\nass
\theta(\s_\infty, V_\infty^+) &=0, \qquad\quad\quad B_\infty^+ = M_2(\R)\\
\nass
\Wald(\s_\infty,\psi_\infty)&=\text{DS}_2=\text{weight 2 disc. series of $GL_2(\R)$.}
\endalign
$$
Since we will be considering only holomorphic cusp forms of weight $\frac32$
and level $4N$ with $N$ odd and square free, these examples give all of the 
relevant local information.

The local root number $\e_p(\frac12,\Wald(\s_p))=\pm1$ and the 
invariant $\e_p(\s_p)$ are related as follows. There is an element
$\bold{-1}$ in the center of $G'_p$ which maps to $-1\in G_p$, and Waldspurger 
defines a sign $\e(\s_p,\psi_p)$, the central sign of $\s_p$, by
$$\s_p(\bold{-1}) = \e(\s_p,\psi_p)\,\chi_\psi(-1)\cdot I_{\s_p},\tag9.11$$
\cite{\waldquaternion}, p.225. Then, 
$$\e_p(\s_p) = \e_p(\frac12,\Wald(\s_p,\psi_p))\,\e(\s_p,\psi_p).\tag9.12$$ 
Note that
$$\prod_{p} \e(\s_p,\psi_p)=1.\tag9.13$$

Thus, for a given global $\s\simeq \tt_p\s_p$, 
$$\align
\e(\frac12,\Wald(\s,\psi))=+1 &\iff 
\text{there is a $B/\Q$ with $\roman{inv}_p(B) =\e_p(\s_p)$ for all
$p$}\tag9.14\\
\nass 
&\iff \text{there is a $B/\Q$ with $\tt_p\theta(\s_p,V^B_p)\ne0$.}
\endalign
$$
The algebra $B$ is then unique, and, if $B'\not\simeq B$, 
the global theta lift $\theta(\s;V^{B'})=0$ for local reasons. On the other hand,
the global theta lift for $V^B$ is 
$$
\theta(\s,V^B) \simeq \cases  \tt_p\theta(\s_p,V^B_p)&\text{or}\\
\nass
0,&{}\endcases\tag9.15
$$
and Waldspurger's beautiful result, \cite{\waldquaternion}, is that
$$\theta(\s,V^B)\ne0 \iff L(\frac12,\Wald(\s))\ne0.\tag9.16
$$

\subheading{\Sec10. The doubling integral}

The key to linking Waldspurger's theory to the arithmetic theta lift is the 
doubling integral
representation of the Hecke L-function, \cite{\psrallis}, \cite{\li}, \cite{\doubling}, 
and in classical language, \cite{\garrettII}, \cite{\boecherer}.
Let $\Gt = \roman{Sp}_2$ be the symplectic group of rank $2$ over $\Q$, and let $\Gt'_\A$ 
be the $2$--fold metaplectic cover of $\Gt(\A)$. Recall that $G=\SL_2 = \roman{Sp}_1$. 
Let $i_0: G\times G\rightarrow\Gt$ be the standard embedding:
$$i_0:\ \pmatrix a_1&b_1\\c_1&d_1\endpmatrix\times\pmatrix a_2&b_2\\c_2&d_2\endpmatrix
\longmapsto \pmatrix a_1&{}&b_1&{}\\{}&a_2&{}&b_2\\c_1&{}&d_1&{}\\{}&c_2&{}&d_2
\endpmatrix.\tag10.1$$
For $g\in G$, let 
$$g^\vee = \roman{Ad}\pmatrix 1&{}\\{}&-1\endpmatrix\cdot g,\tag10.2$$
and let 
$$i(g_1,g_2) = i_0(g_1,g_2^\vee).\tag10.3$$
Let $\Pt\subset \Gt$ be the standard Siegel parabolic, and, for $s\in \C$, let
$\It(s)$ be the induced representation
$$\It(s) = \roman{Ind}_{\Pt'_\A}^{\Gt'_\A}(\delta^{s+\frac32}) \simeq \tt_p \It_p(s),\tag10.4$$
where $\delta$ is a certain lift to the cover $\Pt'_\A$ of the modulus character
of $\Pt(\A)$. Here we are using unnormalized induction. 
For a section $\P(s)\in \It(s)$, there is a
Siegel--Eisenstein  series 
$$E(g',s,\P) = \sum_{\gamma\in\Pt'_\Q\back \Gt'_\Q} \P(\gamma g',s),\tag10.5$$
convergent for $\Re(s)>\frac32$, and with an analytic continuation in $s$ 
satisfying a functional equation relating $s$ and $-s$.  
For $\s \simeq \tt_p\s_p$ a cuspidal representation in $\Cal A_{0}(G')$, as above, the 
doubling integral is defined as follows: \hfb
For vectors $f_1$, $f_2\in \s$ and $\P(s)\in\It(s)$, 
$$Z(s,f_1,f_2,\P) = \int_{G'_\Q\back G'_\A\times G'_\Q\back G'_\A} 
\overline{f_1(g'_1)}\, f_2(g'_2)\,E(i(g_1',g_2'),s,\P)\,dg'_1\,dg'_2.\tag10.6$$  
If $f_1$, $f_2$ and $\P(s)$ are factorizable and unramified outside of 
a set of places $S$, including $\infty$ and $2$, then, \cite{\psrallis}, \cite{\li},
\cite{\krannals}\footnote{Note that \cite{\psrallis} and \cite{\krannals} 
deal with the symplectic case, while the metaplectic case needed here is 
covered in \cite{\li}.},  
$$Z(s,f_1,f_2,\P) = \frac{1}{\zeta^S(2s+2)}\,
L^S(s+\frac12,\Wald(\s))\cdot
\prod_{p\in S} Z_p(s,f_{1,p},f_{2,p},\P_p),\tag10.7$$
where $Z_p(s,f_{1,p},f_{2,p},\P_p)$ is a local zeta integral 
depending on the local components at $p$. 

Now suppose that $\s_\infty = \roman{HDS}_{\frac32}$,  
and take the archimedean local components $f_{1,\infty}= 
f_{2,\infty}$ to be the weight $\frac32$ vectors. We can identify $f_1$ and $f_2$ 
with classical cusp forms of weight $\frac32$. 
Taking $\P_\infty(s)\in \It_\infty(s)$ to be the standard
weight $\frac32$ section $\P_{\infty}^{\frac32}(s)$, we can 
write the Eisenstein series as a classical Siegel 
Eisenstein series $E(\tau,s,\P_f)$ of weight $\frac32$, where $\tau\in\H_2$, the Siegel space 
of genus $2$, 
and $\P_f(s)$ is the finite component of $\P(s)$.  

For a given indefinite quaternion algebra $B$, there is a section $\P_f(s)$ defined as follows. 
For a finite prime $p$, the group $\Gt'_p$ acts on the Schwartz space 
$S((V_p^\pm)^2)$ via the Weil representation $\o$ determined by $\psi_p$, and there 
is a map 
$$\l_p:S((V_p^\pm)^2)\lra \It_p(0), \qquad \ph_p\mapsto \l_p(\ph_p)(g') = \o(g')\ph_p(0).\tag10.8$$
Here $V_p^\pm$ is the ternary quadratic space defined in (5.10).
Note that a section $\P_p(s)$ is determined by its restriction to the 
compact open subgroup $K'_p$, the inverse image of $\Gt(\Z_p)$ in $\Gt'_p$. 
A section is said to be standard if this restriction is independent of $s$.  
The function $\l_p(\ph_p)\in \It_p(0)$ has a unique extension to a standard section 
of $\It_p(s)$. 
Fix a maximal order $R_p^\pm$ in $B_p^\pm$, and let $\ph_p^\pm\in S((V_p^\pm)^2)$ be the 
characteristic function of $(R_p^\pm\cap V_p^\pm)^2$. Also let $R_p^e\subset R_p^+$ 
be the Eichler order\footnote{So $R_p^+ =M_2(\Z_p)$ and 
$x\in R_p^e$ iff $c\equiv 0\mod p$.} of index $p$,
and let
$\ph_p^e$  be the characteristic function of $(R_p^e\cap V_p^+)^2$. 
We then have standard sections $\P_p^0(s)$, $\P_p^-(s)$ and $\P_p^e(s)$ whose restrictions 
to $K'_p$ are $\l_p(\ph_p^+)$, $\l_p(\ph_p^-)$, and $\l_p(\ph_p^e)$ respectively. 
Following \cite{\krinvent}, let
$$\Pht_p(s) = \P_p^-(s) + A_p(s)\P_p^0(s) + B_p(s)\,\P_p^e(s),\tag10.9$$
where $A_p(s)$ and $B_p(s)$ are entire functions of $s$ such that 
$$A_p(0)=B_p(0) =0, \qquad \text{and}\qquad A'_p(0) = -\frac{2}{p^2-1}\,\log(p), \quad
B'_p(0) = \frac12\frac{p+1}{p-1}\,\log(p).\tag10.10$$
Then let
$$\Pht^B(s) = \bigg(\tt_{p\mid D(B)} \Pht_p(s)\bigg) \tt\bigg( \tt_{p\nmid D(B)}
\P_p^0(s)\bigg).\tag10.11$$

Using this section, we can define a normalized Siegel Eisenstein series of weight $\frac32$ and genus 
$2$ attached to $B$ by
$$\Cal E_2(\tau,s;B) = \eta(s,B)\,\zeta_{D(B)}(2s+2)\,E(\tau,s;\Pht^B),\tag10.12$$
where $\eta(s,B)$ is a certain normalizing factor and the partial zeta function 
is as in (2.14). 
Then we have a precise version of the doubling identity, \cite{\doubling}:
\proclaim{Theorem 10.1}\footnote{At the time of this writing, we must assume that 
$2\nmid D(B)$ and that $\e_2(\s_2)=+1$. It should not be 
difficult to remove these restrictions.} 
For every indefinite quaternion algebra $B$ over $\Q$, the associated Siegel--Eisenstein
series 
$\Cal E_2(\tau,s;B)$
of genus $2$ and weight $\frac32$ has the following property. For each holomorphic 
`newform' $f$ of weight $\frac32$ and level $4D(B)$
associated to an irreducible cuspidal representation
$\s=\tt_p\s_p$:
$$\align
&\gs{\,\Cal E_2(\pmatrix
\tau_1&{}\\{}&-\bar\tau_2\endpmatrix,s;B)\,}{\,\overline{f(\tau_2)}\,}_{\text{Pet,
$\tau_2$}}\\
\nass
&\qquad\qquad{}= C(s)\,C(s;\s;B)\,L(s+\frac12,\Wald(\s))\cdot f(\tau_1),
\endalign
$$
where\footnote{These are explicit elementary factors.} $C(0)\ne0$, and
$$C(s;\s;B) = \prod_{p\mid D(B)} C_p(s;\s_p;B_p),$$
with
$$C_p(0;\s_p;B_p) = \cases 1 &\text{if $\e_p(\s_p)=-1$,}\\
\nass
0&\text{if $\e_p(\s_p)=+1$,}
\endcases
$$
and, if $\e_p(\s_p)=+1$, 
$$C'_p(0;\s_p;B_p) = \log(p).$$
\endproclaim

Note that $\e_p(\s_p)=1$ for $p\nmid 4D(B)$, and that
$$\e(\frac12,\Wald(\s)) = -\prod_{p\mid D(B)} \e_p(\s_p).\tag10.13$$
Thus, for example, for $\e(\frac12,\Wald(\s))=+1$, an odd number of $p\mid D(B)$ 
have $\e_p(\s_p)=+1$. 

\proclaim{Corollary 10.2} With the notation and assumptions of
Theorem~10.1, 
$$\align
&\gs{\,\Cal E'_2(\pmatrix
\tau_1&{}\\{}&-\bar\tau_2\endpmatrix,0;B)\,}{\,\overline{f(\tau_2)}\,}_{\text{Pet,
$\tau_2$}}\\
\nass
\nass
&{}= f(\tau_1)\cdot C(0)\cdot  
\cases L'(\frac12,\Wald(\s)),&\text{if $\e_p(\s_p)=-1$}\\ 
{}&\text{for all $p\mid D(B)$,}\\
\nass
L(\frac12,\Wald(\s))&\text{if $\e_p(\s_p)=+1$}\\
{\quad\times\log(p),}&\text{for a unique $p\mid D(B)$,}\\
\nass
0&\text{otherwise.}
\endcases
\endalign
$$
\endproclaim

\subheading{\Sec11. The arithmetic inner product formula} 

Returning to the arithmetic theta lift and arithmetic theta function $\thh(\tau)$, 
there should be a second relation between derivatives of Eisenstein series 
and arithmetic geometry, \cite{\annals}:
\proclaim{Conjecture 11.1}
$$\gs{\thh(\tau_1)}{\thh(\tau_2)} = \Cal E_2'(\pmatrix \tau_1&{}\\{}&-\bar\tau_2\endpmatrix,0;B).$$
\endproclaim
Here recall that $\gs{\ }{\ }$ has been extended to be conjugate linear in the 
second factor.  Additional discussion can be found in \cite{\bourbaki}, \cite{\kudlaICM},
\cite{\kudlaMSRI},
\cite{\kryIV}. This conjecture amounts to identities on Fourier coefficients:
$$\align
&\gs{\Zh(t_1,v_1)}{\Zh(t_2,v_2)}\cdot q_1^{t_1}q_2^{t_2}\tag11.1\\
\nass 
&\qquad\qquad{}=\sum_{\matrix\scr T \in \Sym_2(\Z)^\vee\\ \scr \diag(T) = (t_1,t_2)\endmatrix}
\Cal E'_{2,T}(\pmatrix \tau_1&{}\\{}&\tau_2\endpmatrix,0;B).
\endalign
$$
If $t_1t_2$ is not a square, then any $T\in\Sym_2(\Z)^\vee$ with $\diag(T)=(t_1,t_2)$ 
has $\det(T)\ne0$. Thus, only the nonsingular Fourier coefficients of $\Cal E_2(\tau,s,B)$
contribute to the right hand side of (11.1) in this case. Under the same condition, the 
cycles $\Cal Z(t_1)$ and $\Cal Z(t_2)$ do not meet on the generic fiber, although 
they may have common vertical components in the fibers of bad reduction. 

On the other hand, if $t_1t_2 = m^2$, then the singular matrices 
$T= \pmatrix t_1&\pm m\\ \pm m&t_2\endpmatrix$ 
occur on the right side of (11.1), and the cycles $\Cal Z(t_1)$ and $\Cal Z(t_2)$ 
meet in the generic fiber and have common horizontal components.

The results of \cite{\annals} and \cite{\krinvent} yield the following.

\proclaim{Theorem 11.2} Suppose that $t_1t_2$ is not a square. 
In addition, assume that {\rm(11.11)} {\rm({\it resp.} (11.26))} below holds for $p=2$ if\ 
$2\nmid D(B)$ \ {\rm(}resp. $2\mid D(B)$\,{\rm)}. Then the 
Fourier coefficient identity {\rm(11.1)} holds. 
\endproclaim

We now briefly sketch the proof of Theorem~11.2. The basic idea is that there is a
decomposition of the  height pairing on the left hand side of (11.1) into terms indexed by 
$T\in \Sym_2(\Z)^\vee$ with $\diag(T)=(t_1,t_2)$. One can prove 
identities between terms on the two sides corresponding to a given $T$.

Recalling the modular definition of the 
cycles given in section 3, the intersection $\Cal Z(t_1)\cap \Cal Z(t_2)$ can be viewed as
the locus of triples $(A,\iota,\bold x)$, where $\bold x = [x_1,x_2]$ is a pair of 
special endomorphisms $x_i\in V(A,\iota)$ with $Q(x_i)=t_i$. Associated to $\bold x$ is the 
`fundamental matrix' $Q(\bold x) = \frac12\, ((x_i,x_j))\in \Sym_2(\Z)^\vee$, where 
$(\ ,\ )$ is the bilinear form on the
quadratic lattice $V(A,\iota)$.  Thus we may write
$$\Cal Z(t_1)\cap \Cal Z(t_2) = \coprod_{\matrix\scr T \in \Sym_2(\Z)^\vee\\ \scr \diag(T) =
(t_1,t_2)\endmatrix}
\Cal Z(T),\tag11.2$$ 
where $\Cal Z(T)$ is the locus of triples with $Q(\bold x)=T$. 
Note that the fundamental matrix is always positive semidefinite, since the
quadratic form on $V(A,\iota)$ is positive definite. On the other hand, if 
$\det(T)\ne0$, then $\Cal Z(T)_\Q$ is empty, since the space of special 
endomorphisms $V(A,\iota)$ has rank $0$ or $1$ in characteristic $0$. 

For a nonsingular $T\in \Sym_2(\Q)$, there is a unique global ternary 
quadratic space $V_T$ with discriminant $-1$ which represents $T$; 
the matrix of the quadratic form 
on this space is
$$Q_T = \pmatrix T&{}\\{}&\det(T)^{-1}\endpmatrix.\tag11.3$$
The space $V_T$ is isometric to the space of trace zero elements for some quaternion 
algebra $B_T$ over $\Q$, and the local invariants of this 
algebra must differ from those of the given indefinite $B$ at a finite set of places 
$$\roman{Diff}(T,B) := \{\ p\le\infty\ \mid\ 
\roman{inv}_p(B_T)=-\roman{inv}_p(B)\ \},\tag11.4$$
with $|\roman{Diff}(T,B)|$
even\footnote{This is a slightly different definition than that used 
in \cite{\annals}, where $\infty$ is taken to be in the Diff set 
for $T$ indefinite. Hence the difference in parity. }. 

The nonsingular Fourier coefficients of $\Cal E_2(\tau,s;B)$ have a product 
formula
$$\Cal E_{2,T}(\tau,s;B) = \eta(s;B)\,\zeta_{D(B)}(2s+2)\cdot 
W_{T,\infty}(\tau, s;\frac32)\cdot\prod_p
W_{T,p}(s;\Pht_p^B),\tag11.5$$ 
where $\Pht^B(s) = \tt_p \Pht_p^B(s)$ is given by (10.11). Moreover, for a finite prime $p$,  
$$p\in \roman{Diff}(T,B) \iff W_{T,p}(0,\Pht_p^B)=0,\tag11.6$$
while
$$\ord_{s=0}\, W_{T,\infty}(\tau,s;\frac32) = \cases 0 &\text{if $T>0$, and}\\
\nass
1&\text{if $\sig(T) = (1,1)$ or $(0,2)$.}
\endcases\tag11.7
$$

First suppose that $T>0$, so that $\infty\in \roman{Diff}(T,B)$, since $V^B$ has signature 
$(1,2)$ and hence cannot represent $T$. If $|\roman{Diff}(T,B)|\ge4$, then 
$\Cal Z(T)$ is empty and $\Cal E'_{2,T}(\tau,0;B)=0$, so there is no contribution 
of such $T$'s on either side of (11.1).  If $T>0$ and $|\roman{Diff}(T,B)|=2$, 
then $\roman{Diff}(T,B) = \{\infty, p\}$ for a unique finite prime $p$. In this situation, 
it turns out that $\Cal Z(T)$ is supported in the fiber $\Cal M_p$ at $p$. 
There are two distinct cases:
\roster
\item"{(i)}" If $p\nmid D(B)$, then $\Cal Z(T)$ is a finite set of points in $\Cal M_p$.
\item"{(ii)}" If $p\mid D(B)$, then $\Cal Z(T)$ can be a union of components $Y_p$ of the 
fiber at $p$, with multiplicities. 
\endroster

In case (i), the contribution to the height pairing on the left side of (11.1) 
is $\log(p)$ times the sum of the local multiplicities of points in $\Cal Z(T)$. 
It turns out that all points have the same multiplicity, $e_p(T)$, so that 
the contribution to (11.1) has the form
$$e_p(T)\cdot|\hskip .15pt\Cal Z(T)(\bar\Bbb F_p)\hskip .1pt|.\tag11.8$$
The computation of $e_p(T)$ can be reduced to a special case of a problem 
in the deformation 
theory of $p$--divisible groups which was solved by Gross and Keating, \cite{\grosskeating}. 
As explained in section 14 of \cite{\annals},  
their result yields the formula
$$e_p(T) = \cases
\sum_{j=0}^{\frac{\a-1}2} (\a+\b-4j)\,p^j &\text{ if $\a$ is odd},\\
\nass
\sum_{j=0}^{\frac{\a}2-1} (\a+\b-4j)\,p^j + \frac12(\b-\a+1)\,p^{\frac{\a}2}
&\text{ if $\a$ is even,}
\endcases\tag11.9
$$
where, for $p\ne2$, $T$ is equivalent, via the action of $\GL_2(\Z_p)$, 
to $\diag(\e_1\, p^{\a}, \e_2\,p^\b)$, with $\e_1$, $\e_2\in \Z_p^\times$ 
and $0\le\a\le\b$. The same result holds for $p=2$, but with a slightly 
different definition, \cite{\grosskeating}, of the invariants $\a$ and $\b$ of $T$.  

On the other hand, if $\Diff(T,B)=\{\infty,p\}$ with $p\nmid D(B)$, then 
$$\align
\Cal E'_{2,T}(\tau,0;B)\ &=\ \eta(0;B)\,\zeta_{D(B)}(2)\cdot
W_{T,\infty}(\tau,0;\frac32)\tag11.10\\
\nass
{}&\qquad\qquad 
\times \frac{W'_{T,p}(0,\P^0_p)}{W_{T,p}(0,\P_p^-)}\cdot 
\bigg(\, W_{T,p}(0,\P_p^-)\cdot\prod_{\ell\ne p} W_{T,\ell}(0,\Pht_\ell^B)\,\bigg).
\endalign
$$
Here recall that the local sections $\Pht_\ell^B(s)$ are as defined in (10.9) and (10.11).
For $p\ne2$, formulas of Kitaoka, \cite{\kitaoka}, for representation densities 
$\a_p(S,T)$ of binary forms $T$ by unimodular quadratic forms $S$ can be used 
to compute the derivatives of the local Whittaker functions, and yield,
\cite{\annals}, Proposition~8.1 and Proposition~14.6,  
$$\frac{W'_{T,p}(0,\P^0_p)}{W_{T,p}(0,\P_p^-)} = \frac12(p-1)\,\log(p)\cdot e_p(T),\tag11.11$$
where $e_p(T)$ is precisely the multiplicity (11.9)!
On the other hand, up to simple constants\footnote{Hence the notation $\circeq$.},  
$$\align
W_{T,\infty}(\pmatrix \tau_1&{}\\{}&\tau_2\endpmatrix,0;\frac32)\ &\circeq\ 
q_1^{t_1}q_2^{t_2},\tag11.12\\
\noalign{and}
\nass
\nass
W_{T,p}(0,\P_p^-)\cdot\prod_{\ell\ne p} W_{T,\ell}(0,\Pht_\ell^B)\  &\circeq\ 
|\hskip .15pt\Cal Z(T)(\bar\Bbb F_p)\hskip .1pt|.\tag11.13
\endalign
$$
Note that to obtain (11.11) in the case $p=2$, one needs to extend Kitaoka's 
representation density formula to this case; work on this is in progress. 

Next, we turn to case (ii), where the component $\ZZ(T)$ 
is attached to $T$ with $\diag(T)=(t_1,t_2)$ 
and $\Diff(T,B)=\{\infty,p\}$ with $p\mid D(B)$. 
This case is studied in detail in \cite{\krinvent} under the
assumption that
$p\ne2$,  using the $p$--adic uniformization described in section 4 above. 
First, we can base change to $\Z_{(p)}$ and
use the intersection theory explained in section 4
of \cite{\krinvent}. The contribution of $\ZZ(T)$ to the height pairing is then
$$
\chi(\,\ZZ(T),\Cal O_{\ZZ(t_1)}\overset{\Bbb L}\to{\tt}\Cal O_{\ZZ(t_2)}\,)\cdot \log(p).
\tag11.14
$$
Here $\chi$ is the Euler--Poincar\'e characteristic and, 
for quasicoherent sheaves $\Cal F$ and $\Cal G$ on 
$\Cal M\times_{\Spec(\Z)}\Spec(\Z_{(p)})$, with $\roman{supp}(\Cal F)\cap
\roman{supp}(\Cal G)$ contained in the special fiber and proper over $\Spec(\Z_{(p)})$, 
$$\chi(\Cal F\overset{\Bbb L}\to{\tt}\Cal G)\  =\ 
\chi(\Cal F\tt\Cal G)-\chi(\Cal T\!or_1(\Cal F,\Cal G)) + 
\chi(\Cal T\!or_2(\Cal F,\Cal G)).\tag11.15$$

Recall that, for $i=1$, $2$,  $\widehat{\Cal C}_p(t_i)$ 
is the base change to $W$ of the formal 
completion of $\ZZ(t_i)$ along its fiber at $p$. Similarly, we write
$\widehat{\Cal C}_p(T)$ for the analogous formal scheme over $W$ determined 
by $\ZZ(T)$. By Lemma~8.4 of \cite{\krinvent}, 
$$
\chi(\,\ZZ(T),\Cal O_{\ZZ(t_1)}\overset{\Bbb L}\to{\tt}\Cal O_{\ZZ(t_2)}\,)\ =\ 
\chi(\,\widehat{\Cal C}_p(T),\Cal O_{\widehat{\Cal C}_p(t_1)}\overset{\Bbb L}\to{\tt}
\Cal O_{\widehat{\Cal C}_p(t_2)}\,),\tag11.16$$
so that we can calculate after passing to the formal situation. The same arguments
which yield the $p$--adic uniformization, Proposition~4.2, 
of $\widehat{\Cal C}_p(t)$ yields a diagram
$$\matrix
\widehat{\Cal C}_p(T) & \isoarrow & H'(\Q)\back (\star\star)\\
\nass
\downarrow&{}&\downarrow\\
\nass 
\widehat{\Cal M}_p & \isoarrow & H'(\Q)\back \bigg(\ \Cal D^\bullet\times H(\A_f^p)/K^p\ \bigg)
\endmatrix\tag11.17
$$
where
$$(\star\star):=\left\{\ (\bold y,(X,\rho),gK^p)\ \bigg\vert\ \matrix (i)\quad Q(\bold y)=T\\
\nass
 (ii)\quad(X,\rho)\in Z^\bullet(j(\bold y))\\
\nass 
(iii)\quad \bold y\in \big(\,g\,(V(\A_f^p)\cap \hat{O}_B^p)\,g^{-1}\big)^2\endmatrix\
\right\}.\tag11.18$$
Proceeding as in Remark 4.3, we obtain
$$\widehat{\Cal C}_p(T) \isoarrow \big[\,\Gamma'\back \Cal D_T^\bullet\,\big],\tag11.19$$
where
$$\Cal D_T^\bullet = \coprod_{\matrix\scr \bold y\in (L')^2\\ \scr Q(\bold y)=T 
\endmatrix} Z^\bullet(j(\bold y)).\tag11.20$$
Here, recall from (4.25) that 
$$\Gamma' = H'(\Q)\cap \big(\,H'(\Q_p)\times K^p\,) =
\big(\,O_B'\big[\frac{1}{p}\big]\,\big)^\times.\tag11.21$$
Since any $\bold y\in L'$ with $Q(\bold y)=T$ spans a nondegenerate $2$--plane
in $V'$, 
$$\Gamma'_{\bold y} = \Gamma'\cap Z(\Q) \simeq \big(\,\Z\big[\frac{1}{p}\big]\,\big)^\times
=\{\pm1\}\times p^{\Z},\tag11.22$$
and the central element $p$ acts on $\Cal D^\bullet$ by translation by $2$, 
i.e. carries $\Cal D^i$ to $\Cal D^{i+2}$. Thus, unfolding, as in \cite{\krinvent}, p.216, 
we have
$$\chi(\,\widehat{\Cal C}_p(T),\Cal O_{\widehat{\Cal C}_p(t_1)}\overset{\Bbb L}\to{\tt}
\Cal O_{\widehat{\Cal C}_p(t_2)}\,) = 
\sum_{\matrix \scr \bold y\in (L')^2\\
\scr Q(\bold y)=T\\ \scr \mod \Gamma'\endmatrix} 
\chi(\,Z(\bold j),\Cal O_{Z(j_1)}\overset{\Bbb L}\to{\tt}
\Cal O_{Z(j_2)}\,).\tag11.23
$$
Here we have used the orbifold convention, which introduces a factor of $\frac12$
from the $\pm1$ in $\Gamma'_{\bold y}$ which acts trivially, and the fact that the 
two `sheets' $Z(\bold j)=Z^0(\bold j)$ and $Z^1(\bold j)$ make the same contribution. 

Two of the main results of \cite{\krinvent}, Theorem~5.1 and~6.1, give the
following:
\proclaim{Theorem 11.3} (i) The quantity
$$e_p(T):=\chi(\,Z(\bold j),\Cal O_{Z(j_1)}\overset{\Bbb L}\to{\tt}
\Cal O_{Z(j_2)}\,)$$
is the intersection number, \cite{\krinvent}, section 4,
of the cycles $Z(j_1)$ and $Z(j_2)$ in
the  formal scheme $\Cal D_p$. It
depends only on the $\GL_2(\Z_p)$--equivalence class of $T$. \hfb
(ii) For $p\ne 2$, and for $T\in \Sym_2(\Z_p)$ which is $\GL_2(\Z_p)$--equivalent to
$\diag(\e_1 p^\a,\e_2 p^\b)$, with $0\le\a\le\b$ and $\e_1$, $\e_2\in \Z_p^\times$,  
$$e_p(T) = \a+\b+1-\cases
p^{\a/2}+2\frac{p^{\a/2}-1}{p-1}&\text{if $\a$ is even and $(-\e_1,p)_p=-1$,}\\
\nass
(\b-\a+1)\,p^{\a/2} + 2\frac{p^{\a/2}-1}{p-1}&\text{if $\a$ is even and $(-\e_1,p)_p=1$,}\\
\nass
2\frac{p^{(\a+1)/2}-1}{p-1}&\text{if $\a$ is odd.}
\endcases
$$
\endproclaim

Part (i) of Theorem~11.3, together with (11.16) and (11.23), yields
$$\chi(\,\ZZ(T),\Cal O_{\ZZ(t_1)}\overset{\Bbb L}\to{\tt}\Cal O_{\ZZ(t_2)}\,)\ 
=\ 
e_p(T)\cdot \big(\,\sum_{\matrix \scr \bold y\in (L')^2\\
\scr Q(\bold y)=T\\ \scr \mod \Gamma'\endmatrix} 1\,\big),\tag11.24$$
which is the analogue of 
(11.8) in the present case.

On the other hand, if $\Diff(T,B)=\{\infty,p\}$ with $p\mid D(B)$, the term on the 
right side of (11.1) is
$$\Cal E'_{2,T}(\tau,0;B) =c\cdot W_{T,\infty}(\tau,0;\frac32) 
\cdot W'_{T,p}(0,\Pht_p)\cdot \prod_{\ell\ne p} W_{T,\ell}(0,\Pht_\ell^B).
\tag11.25$$
where $c=\eta(0;B)\,\zeta_{D(B)}(2)$.
By \cite{\krinvent}, Corollary~11.4, the section $\Pht_p(s)$ in (10.9) satisfies,  
$$\align
W_{T,p}'(0,\Pht_p) &= p^{-2}(p+1)\,\log(p)\cdot e_p(T)\tag11.26\\
\noalign{while}
\nass
\prod_{\ell\ne p} W_{T,\ell}(0,\P_\ell) &\circeq 
\big(\,\sum_{\matrix \scr \bold y\in (L')^2\\
\scr Q(\bold y)=T\\ \scr \mod \Gamma'\endmatrix} 1\,\big).\tag11.27
\endalign
$$
can be thought of as the number of `connected components' of $\ZZ(T)$. 
Note that, the particular choice (10.10) of the coefficients $A_p(s)$ and $B_p(s)$ 
in the definition of $\Pht_p(s)$ was dictated by the identity (11.16), the proof of
which is based on results of Tonghai Yang, \cite{\yangden}, 
on representation densities $\a_p(S,T)$ of 
binary forms $T$ by nonunimodular forms $S$. 
To obtain (11.16) in the case $p=2$, one needs to extend both the intersection 
calculations of \cite{\krinvent} and the density formulas of \cite{\yangden} to this case. 
The first of these tasks, begun in the appendix to section 11 of \cite{\kryIII}, is now
complete. The second is in progress. 

Finally, the contribution to the right side of (11.1) 
of the terms for $T$ of signature $(1,1)$ or $(0,2)$, 
which can be calculated using the formulas of \cite{\shimuraconf}, coincides with the
contribution of the star product of the Green functions $\Xi(t_1,v_1)$ 
and $\Xi(t_2,v_2)$ to the height pairing. This is a main result of \cite{\annals}; 
for a sketch of the ideas involved, cf. \cite{\bourbaki}. 

This completes the sketch of the proof of Theorem~11.2. 
\qed\!\!\!\qed

Before turning to consequences, we briefly sketch how part (ii) of 
Theorem~11.3 is obtained. By part (i) of that Theorem, it will suffice to 
compute the intersection number $(Z(j_1),Z(j_2))$ for $j_1$ and $j_2\in V'(\Q_p)$ 
with $j_1^2=-Q(j_1)=-\e_1p^\a$, $j_2^2=-Q(j_2)=-\e_2 p^\b$ and 
$(j_1,j_2) = j_1j_2+j_2j_1=0$. By Theorem~5.1 of \cite{\krinvent}, 
$(Z(j_1),Z(j_2)) = (Z(j_1)^{\roman{pure}},Z(j_2)^{\roman{pure}}),$ 
so that we may use the description of $Z(j)^{\roman{pure}}$ given in 
Proposition~4.5, above. 

The following result, which combines Lemmas~4.7, 4.8, and 4.9 of 
\cite{\krinvent}, describes the intersection numbers of individual components.
Recall that the vertical components $\Bbb P_{[\L]}$ of $\Cal D$ are indexed 
by vertices $[\L]$ in the building $\Cal B$ of $\roman{PGL}_2(\Q_p)$. 
\proclaim{Lemma 11.4} (i) For a pair of vertices $[\L]$ and $[\L']$, 
$$\align
(\Bbb P_{[\L]},\Bbb P_{[\L']}) &= \cases 1&\text{if $([\L],[\L'])$ is an edge,}\\
\nass
-(p+1)&\text{if $[\L]=[\L']$, and}\\
\nass
0&\text{otherwise.}
\endcases\\
\noalign{(ii)}
(Z(j_1)^h,Z(j_2)^h) &= \cases 1&\text{if $\a$ and $\b$ are odd,}\\
0&\text{otherwise.}
\endcases\\
\noalign{(iii)}
(Z(j_1)^h,\Bbb P_{[\L]}) &= \cases 
2&\text{if $\a$ is even, $(-\e_1,p)_p=-1$, and $\Cal B^{j_1}={[\L]}$,}\\
\nass
1&\text{if $\a$ is odd and $d([\L],\Cal B^{j_1})=\frac12$, and}\\
\nass
0&\text{otherwise.}
\endcases
\endalign
$$
and similarly for $Z(j_2)^h$. 
\endproclaim
The computation of the intersection number $(Z(j_1)^{\roman{pure}},Z(j_2)^{\roman{pure}})$
is thus reduced to a combinatorial problem. Recall from Proposition~4.5 
above that the multiplicity in
$Z(j_1)^{\roman{pure}}$ of a vertical component $\Bbb P_{[\L]}$
indexed by a vertex
$[\L]$ is determined by the distance of $[\L]$ from the fixed point set $\Cal B^{j_1}$ 
of $j_1$ on $\Cal B$ by the formula
$$\mu_{[\L]}(j_1) = \max\{\,0,\,\frac{\a}2 - d([\L],\Cal B^{j_1})\,\}.$$
In particular, this multiplicity is zero outside of the tube $\Cal T(j_1)$ 
of radius $\frac{\a}2$ around $\Cal B^{j_1}$. Of course, the 
analogous description holds for $Z(j_2)^{\roman{pure}}$.
Our assumption that the matrix 
of inner products of $j_1$ and $j_2$ is diagonal implies that $j_1$ and $j_2$ 
anticommute, and hence the relative position of the fixed point sets 
and tubes $\Cal T(j_1)$ and $\Cal T(j_2)$ is particularly 
convenient.  

For example, consider the case in which $\a$ and $\b$ 
are both even, with $(-\e_1,p)_p=1$ and $(-\e_2,p)_p=-1$. These conditions
mean that $\Q_p(j_1)^\times$ is a split Cartan in $\GL_2(\Q_p)$ and $\Q_p(j_2)^\times$ 
is a nonsplit, unramified, Cartan. The fixed point set $\Cal A=\Cal B^{j_1}$ 
is an apartment, and $\Cal T(j_1)$ is a tube of radius $\frac{\a}2$ 
around it. Moreover, $Z(j_1)^h=\emptyset$, so that $Z(j_1)$ 
consists entirely of vertical components. The fixed point
set
$\Cal B^{j_2}$ is a vertex $[\L_0]$ and $\Cal T(j_2)$ 
is a ball of radius $\frac{\b}2$ around it. Since $j_1$
and $j_2$ anticommute, the vertex $[\L_0]$ lies in the apartment $\Cal A=\Cal B^{j_1}$. 
For any vertex $[\L]$, the geodesic from $[\L_0]$ to $[\L]$ 
runs a distance $\ell$ inside the apartment $\Cal A$ and then a distance $r$ 
outside of it. By (i) of Lemma~11.4, the contribution to  
the intersection number of the vertices with $\ell=0$  is
$$-(p+1)\, \frac{\a}2 + (1-p)\sum_{r=1}^{\a/2-1} (\frac{\a}2-r)(p-1)p^{r-1}
\ =\ 1-\a-p^{\a/2}.\tag11.28$$
Here the first term is the contribution of $[\L_0]$, since
$$(\Bbb P_{[\L_0]},Z(j_2)^v) = -(p+1)\frac{\b}2 + (p+1)(\frac{\b}2-1) = -(p+1),\tag11.29$$
where $Z(j_2)^v$ is the vertical part of $Z(j_2)^{\roman{pure}}$. 
Similarly, for $r>0$, there are $(p-1)p^{r-1}$ vertices $[\L]$ at distance $r$ (with $\ell=0$), 
and each contributes
$$(\Bbb P_{[\L]},Z(j_2)^v) = (\frac{\b}2-r+1) -(p+1)(\frac{\b}2-r) + p(\frac{\b}2-r-1)\ 
=\ 1-p.\tag11.30$$ 
Here we have used the fact that all such vertices lie inside the ball
$\Cal T(j_2)$. Next, consider
vertices with
$1\le\ell\le(\b-\a)/2$. Note that,  if such a
vertex $[\L]$ has
$r=\frac{\a}2$, then $[\L]\in\Cal T(j_2)$. The contribution of vertices with
$\ell$ values in this range is
$$2(1-p)\sum_{\ell=1}^{(\b-\a)/2} \bigg(\,\frac{\a}2+\sum_{r=1}^{\a/2-1}
(\frac{\a}2-r)(p-1)p^{r-1}\,\bigg)  =(\a-\b)(p^{\a/2-1}-1).\tag11.31$$
Next, in the range $(\b-\a)/2<\ell<\b/2$ and with $r<\b/2-\ell-1$, $[\L]\in\Cal T(j_2)$,
and  such $[\L]$'s contribute
$$2(1-p)\sum_{\ell=(\b-\a)/2+1}^{\b/2-1} \bigg(\,\frac{\a}2+\sum_{r=1}^{\b/2-\ell-1}
(\frac{\a}2-r)(p-1)p^{r-1}\,\bigg)  =2\a-4\,\frac{p^{\a/2-1}-1}{p-1}.\tag11.32$$
Finally, the vertices on the boundary of the ball $\Cal T(j_2)$ contribute
$$\a+2\sum_{\ell=(\b-\a)/2+1}^{\b/2-1} (\frac{\a}2-(\frac{\b}2-\ell))(p-1)p^{\b/2-\ell-1}
\ =\ 2\,\frac{p^{\a/2-1}-1}{p-1}.\tag11.33$$
Summing these contributions and adding $(\Bbb P_{[\L_0]},Z(j_2)^h)\cdot\a/2 = \a$, 
we obtain the
quantity claimed in (ii) of Theorem~11.3  in this case. The other cases are similar. \qed

We now describe the consequences of Theorem~11.2. 
\proclaim{Corollary 11.5} Assume that the $p$--adic density 
identity (11.11) ( resp. (11.26)) holds for $p=2$ when $2\mid D(B)$ (resp. $2\nmid D(B)$).
Then
$$\align
\gs{\thh(\tau_1)}{\thh(\tau_2)} &= \ \Cal E_2'(\pmatrix \tau_1&{}\\{}&-\bar\tau_2\endpmatrix,0;B)\\
\nass
{}&\qquad + \sum_{t_1}\ \bigg(\,\sum_{\matrix \scr t_2\\ \scr t_1t_2=\roman{square}\endmatrix} 
c(t_1,t_2,v_1,v_2)\cdot\bar q_2^{t_2}\, \bigg)\,q_1^{t_1}.
\endalign
$$
for some coefficients $c(t_1,t_2,v_1,v_2)$. 
\endproclaim

The extra term on the right hand side should vanish, according to Conjecture~11.1, 
but, in any case, the coefficient of each $q_1^{t_1}$ is a modular form 
of weight $\frac32$ in $\tau_2$ 
with only 
one square class of Fourier coefficients, i.e., is a distinguished form. But such forms
all come from $O(1)$'s, and so are orthogonal to cusp forms in $\Cal A_{00}(G')$, 
\cite{\gelbartps}. 
\proclaim{Corollary 11.6} With the notation and assumptions of Theorem~10.1
and assuming that the $p$--adic density identity (11.11) 
holds for $p=2$, 
$$\gs{\thh(\tau_1)}{\thh(f)} = f(\tau_1)\cdot C(0)\cdot  
\cases L'(\frac12,\Wald(\s)),&\\ 
\nass
L(\frac12,\Wald(\s))\cdot \log(p),&\\
\nass
0&
\endcases.
$$
\endproclaim
\demo{Proof}
$$\align
\gs{\thh(\tau_1)}{\thh(f)} &= \gs{\thh(\tau_1)}{\gs{\thh}{f}_{\roman{Pet}}}\\
\nass
{}&=\big\langle\,\gs{\thh(\tau_1)}{\thh(\tau_2)},
\overline{f(\tau_2)}\ \big\rangle_{\roman{Pet}}\tag11.34\\
\nass
{}&= \big\langle\,\Cal E_2'(\pmatrix
\tau_1&{}\\{}&-\bar\tau_2\endpmatrix,0;B), \overline{f(\tau_2)}\ \big\rangle_{\roman{Pet}}\\
\nass
&{}= f(\tau_1)\cdot C(0)\cdot  
\cases L'(\frac12,\Wald(\s)),&\\ 
\nass
L(\frac12,\Wald(\s))\cdot \log(p),&\\
\nass
0&
\endcases
\endalign
$$
\qed\enddemo

\proclaim{Corollary 11.7} With the notation and assumptions 
of Corollary~11.6, \hfb
(i) If $\e_p(\s_p)=+1$ for more than
one $p\mid D(B)$, 
then $\thh(f)=0$.\hfb
Otherwise\hfb
(ii) If $\e(\frac12,\Wald(\s))=+1$, then there is a unique 
prime $p\mid D(B)$ such that\hfb
$\e_p(\s_p)=+1$, $\thh(f)\in \Ver_p \subset \Ver$, and 
$$\gs{\thh(f)}{\thh(f)} = C(0)\cdot \gs{f}{f}\cdot L(\frac12,\Wald(\s))\cdot\log(p).$$
(iii) If $\e(\frac12,\Wald(\s))=-1$, then $\thh(f)\in \MWt$ and 
$$\gs{\thh(f)}{\thh(f)} = C(0)\cdot \gs{f}{f}\cdot L'(\frac12,\Wald(\s)).$$
\endproclaim

This result gives an analogue of Waldspurger's theory, 
described in section 9 above, and of the Rallis inner product
formula, 
\cite{\rallisinnerprod}, for the 
arithmetic theta lift (8.1). Of course, the 
result only applies to forms of weight $\frac32$, and 
we have made quite strong restrictions on the level of $f$. 
The restrictions 
on the level can most likely be removed with some additional work, but the 
restriction on the weight is essential to our setup.  

In the case $\e(\frac12,\Wald(\s))=+1$, part (ii) gives a  geometric
interpretation of the value
$L(\frac12,\Wald(\s))$ 
in terms of vertical components of $\Cal M$, analogous to the geometric
interpretation of the value of the base change L-function given by Gross in his Montreal paper,
\cite{\grossmont}.  

Finally, if $\thh(f)\ne0$, let
$$\Zh(t)(f)=\Zh(t,v)(f) = \frac{\gs{\Zh(t,v)}{\thh(f)}}{\gs{\thh(f)}{\thh(f)}}\cdot \thh(f)
\tag11.35$$
be the component of the cycle $\Zh(t,v)$ along the line spanned by $\thh(f)$. 
Note that, by Corollary~11.6, this projection does not depend on $v$. 
 The following is an analogue of the Gross-Kohnen-Zagier relation, \cite{\grosskohnenzagier}, 
\cite{\zagier}. 
Note that it holds in both cases
$\e(\frac12,\Wald(\s))=+1$ or $\e(\frac12,\Wald(\s))=-1$. 

\proclaim{Corollary 11.8}
$$\sum_{t} \Zh(t)(f)\cdot q^t = \frac{f(\tau)\cdot \thh(f)}{\gs{f}{f}}.$$
\endproclaim
Thus, the Fourier coefficients of $f$ encode the position of 
the cycles $\Zh(t)(f)$ on the $\thh(f)$ line.
Both sides are invariant under scaling of $f$ by a constant factor.

For the related results of Gross--Kohnen--Zagier in the case $\e(\frac12,\Wald(\s))=-1$, cf.
\cite{\grosskohnenzagier}, especially Theorem C, p.~503, the discussion on
pp.~556--561, and the 
examples in \cite{\zagier}, where the analogies with the work of 
Hirzebruch--Zagier are also explained and used in the proof! Our class
$\thh(f)/<f,f>$, which arises as the image of $f$ under the arithmetic theta lift, is the
analogue of the class $y_f$ in \cite{\grosskohnenzagier}. To proved the main part of 
Theorem C of \cite{\grosskohnenzagier} in our case, namely that the $\pi_f$--isotypic
components, where
$\pi=\Wald(\s)$, of the  classes $\Zh(t,v)$ either vanish or lie on a line, we can use the Howe
duality Theorem 
\cite{\howe}, \cite{\waldspurgerLHD} for the
local theta correspondence, cf. \cite{\ariththeta}. 

In \cite{\grossMSRI}, Gross gives a beautiful representation theoretic
framework in which the
Gross--Zagier formula, \cite{\grosszagier}, and a result of Waldspurger, 
\cite{\waldcentral}, can be viewed together. He works with 
unitary similitude groups $G=G^B=GU(2)$, constructed from the choice of a quaternion
algebra $B$ over a global field $k$ and a  
quadratic extension $E/k$ which splits $B$. The torus $T$ over $k$ with 
$T(k)=E^\times$ embeds in 
$$G^B(k) = (B^\times\times E^\times)/\Delta k^\times.$$
For a place $v$ of $k$, the results of Tunnel, \cite{\tunnell}, and Saito, \cite{\saito}, 
show that 
the existence of local $T(k_v)$--invariant functionals on 
an irreducible admissible representation $\Pi_v$ of $G_v$ is
controlled by the local root number $\e_v(\Pi_v)$ of the Langlands L-function 
defined by a $4$ dimensional symplectic representation of ${}^LG$,
\cite{\grossMSRI}, section 10.  
There is a dichotomy phenomenon. If $v$ is a place which is not
split in 
$E/k$, then the local quaternion algebras $B_v^+$ and $B_v^-$ 
are both split by $E_v$, so that the torus $T_v$ embeds in both 
similitude groups $G_v^+$ 
and $G_v^-$. If $\Pi_v$ is a
discrete series representation of $G_v^+$, then, by the local 
Jacquet--Langlands correspondence, there is an associated representation 
$\Pi_v'$ of $G_v^-$, and
$$\dim \roman{Hom}_{T_v}(\Pi_v,\C) + \dim \roman{Hom}_{T_v}(\Pi_v',\C) =1.$$
For a global cuspidal automorphic representation $\Pi$ of 
$G(\A) = (\GL_2(\A)\times E^\times_\A)/\A^\times$, there is a finite 
collection of quaternion algebras $B$ over $k$, 
split by $E$, with automorphic cuspidal 
representations $\Pi^B$ associated to $\Pi$ by the global 
Jacquet--Langlands correspondence.  When the global
root number
$\e(\Pi)$ for the degree $4$ Langlands L-function $L(s,\Pi)$ 
is $+1$, Waldspurger's  theorem,
\cite{\waldcentral}, says that the nonvanishing  of 
the automorphic $T(\A)$--invariant linear functional, defined
on $\Pi^B$ 
by integration over 
$\A^\times T(k)\back T(\A)$, is equivalent to (i) the nonvanishing of the 
local linear functionals\footnote{By local root number 
conditions, this uniquely determines the quaternion 
algebra $B$ for the given $E$ and $\Pi$.} and (ii) the 
nonvanishing of the central value $L(\frac12,\Pi)$ of the L-function.  
Suppose that $k$ is a totally real number field, $E$ is a
totally imaginary quadratic extension, and $\Pi_v$ is a discrete series 
of weight $2$ at every archimedean place of $k$. Then, 
in the case $\e(\Pi)=-1$, the local root number conditions
determine a Shimura curve $M^B$ over $k$, and Gross defines a $T(\A_f)$--invariant 
linear functional on the associated Mordell--Weil space $\roman{MW}^B$
using the $0$--cycle attached to $T$. Gross conjectures that
the nonvanishing of this functional on the $\Pi^B_f$ component of $\roman{MW}^B$
is equivalent to the nonvanishing of the central derivative $L'(\frac12,\Pi)$, 
and that there is an explicit expression for this quantity in terms 
of the height pairing of a suitable `test vector'. 
In the case $k=\Q$, this is the classical Gross--Zagier formula, \cite{\grosszagier}, and 
in there is work of Zhang in the general case \cite{\zhang}. 

Thus, there is a close parallel between Gross's `arithmetic' version 
of Waldspurger's central value result \cite{\waldcentral} and our 
`arithmetic' version of Waldspurger's results on the Shimura lift.  
It should be possible to formulate our conjectures 
above for Shimura curves over an arbitrary totally real field $k$. 
It would be interesting to find a direct connection between our constructions
and the methods used by Shou-Wu Zhang, \cite{\zhang}. 

Finally, it is possible to formulate a similar theory for central value/derivative of the 
triple product L-function. This is discussed in \cite{\grosskudla}, 
cf. also \cite{\harriskudlaII}, \cite{\harriskudlaIII}. 
It may be that all of these examples can be covered by some 
sort of arithmetic version of Jacquet's relative trace formula, 
\cite{\jacquet}, \cite{\jacquetwald}, \cite{\baruchmao}. 
Further speculations about `arithmetic theta functions' and derivatives 
of Eisenstein series can be found in \cite{\kudlaICM} and \cite{\kudlaMSRI}.

\vskip.5in

\redefine\vol{\oldvol}


\Refs 
\widestnumber\key{44}
\parskip=10pt

\ref\key{\baruchmao}
\by E. M. Baruch and Z. Mao
\paper Central value of automorphic L-functions
\jour preprint
\yr 2002
\endref

\ref\key{\boecherer}
\by S. B\"ocherer
\paper \"Uber die Funktionalgleichung automorpher L-Funktionen zur
Siegelschen Modul-gruppe
\jour J. reine angew. Math. 
\vol 362 
\yr 1985
\pages 146--168
\endref 


\ref\key{\borchinventI}
\by R. Borcherds
\paper  Automorphic forms on $\text{\rm O}_{s+2,2}(\text{\bf R})$ and infinite products
\jour Invent. math.
\vol 120
\yr 1995
\pages 161--213
\endref

\ref\key{\borchinventII}
\bysame
\paper  Automorphic forms with singularities on Grassmannians
\jour Invent. math.
\vol 132
\yr 1998
\pages 491--562
\endref


\ref\key{\borchduke}
\bysame
\paper The Gross-Kohnen-Zagier theorem in higher dimensions
\jour Duke Math. J.
\yr 1999
\vol 97
\pages 219--233
\endref

\ref\key{\borchdukeII}
\bysame
\paper Correction to: ``The Gross-Kohnen-Zagier theorem in higher dimensions''
\jour Duke Math. J.
\yr 2000
\vol 105
\pages 183--184
\endref


\ref\key{\bost}
\by J.--B. Bost
\paper Potential theory and Lefschetz theorems for arithmetic surfaces
\jour Ann. Sci. \'Ecole Norm. Sup.
\yr 1999
\vol 32
\pages 241--312
\endref

\ref\key{\bostUMD}
\bysame
\paper\rm Lecture at Univ. of Maryland, Nov. 11, 1998
\endref

\ref\key{\boutotcarayol}
\by J.-F. Boutot and H. Carayol
\paper Uniformisation p-adique des courbes de Shimura
\inbook Courbes Modulaires et Courbes de Shimura
\bookinfo Ast\'erisque, vol. {\bf196--197}
\yr 1991
\pages 45--158
\endref

\ref\key{\bruinierI}
\by J. H.  Bruinier
\paper  Borcherds products and Chern classes of Hirzebruch--Zagier divisors
\jour Invent. Math. 
\yr 1999
\vol138
\pages 51--83
\endref

\ref\key{\bruinierII}
\bysame
\book  Borcherds products on O(2,l) and Chern classes of Heegner divisors
\bookinfo Lecture Notes in Math. 1780
\yr 2002
\publ Springer
\publaddr New York
\endref


\ref\key{\brkuehn}
\by J. H. Bruinier, J. I. Burgos Gil, and U. K\"uhn
\paper Borcherds products in the arithmetic intersection theory of 
Hilbert modular surfaces
\jour in preparation
\endref

\ref\key{\bkk}
\by J. I. Burgos Gil, J. Kramer and U. K\"uhn
\paper  Cohomological arithmetic Chow rings
\jour preprint 
\yr 2003
\vol
\pages
\endref


\ref\key{\cohen}
\by H. Cohen
\paper Sums involving the values at negative integers of L-functions of quadratic characters
\jour Math. Ann. 
\yr 1975
\vol 217
\pages 271--285
\endref

\ref\key{\faltings}
\by G. Faltings
\paper Calculus on arithmetic surfaces
\jour Annals of Math. 
\vol 119
\yr 1984
\pages 387--424
\endref


\ref\key{\funkethesis} 
\by J. Funke  
\paper Rational quadratic divisors and automorphic forms
\jour Thesis, University of Maryland
\yr 1999
\endref

\ref\key{\funkecompo} 
\bysame  
\paper Heegner Divisors and non-holomorphic modular forms
\jour Compositio Math. 
\yr to appear
\vol
\pages
\endref

\ref\key{\garrettII}  
\by P. Garrett
\paper Pullbacks of Eisenstein series; applications
\inbook Automorphic Forms of Several Variables
\bookinfo Tanaguchi Symposium, Katata, 1983
\publ Birkh\"auser
\publaddr Boston
\yr 1984
\endref

\ref\key{\gelbartps}
\by S. Gelbart and I. I. Piatetski-Shapiro
\paper Distinguished representations and modular forms of 
half-integral weight
\jour Invent. Math. 
\yr 1980
\vol 59
\pages 145--188
\endref



\ref\key{\gsihes}
\by H. Gillet and C. Soul\'e
\paper Arithmetic intersection theory
\jour Publ. Math. IHES
\yr 1990
\vol 72
\pages 93--174
\endref


\ref\key{\grossmont}
\by B. H. Gross
\paper Heights and special values of L-series
\inbook Number Theory (Montreal, 1985)
\ed H. Kisilevsky and J. Labute
\bookinfo CMS Conf, Proc. 7
\publ AMS
\publaddr Providence
\yr 1987
\pages 115--187
\endref

\ref\key{\grossMSRI}
\bysame
\paper Heegner points and representation theory
\jour Proc. of Conference on Special Values of Rankin L-Series, MSRI, 
Dec. 2001, to appear.
\endref

\ref\key{\grosskeating}
\by B. H. Gross and K. Keating
\paper On the intersection of modular correspondences
\jour Invent. Math.
\vol 112
\yr 1993
\pages 225--245
\endref

\ref\key{\grosskohnenzagier}
\by B. H. Gross, W. Kohnen and D. Zagier
\paper Heegner points and derivatives of L-functions. II
\jour Math. Annalen
\vol 278
\yr 1987
\pages 497--562
\endref

\ref\key{\grosskudla}
\by B. H. Gross and S. Kudla 
\paper Heights and the central critical values of triple product $L$-functions 
\jour Compositio Math. 
\vol 81
\yr 1992
\pages 143--209
\endref


\ref\key{\grosszagier}
\by B. H. Gross and D. Zagier 
\paper Heegner points and the derivatives of $L$-series
\jour Inventiones math. 
\vol 84
\yr 1986
\pages 225--320
\endref

\ref\key{\harriskudlaII} 
\by M. Harris and S. Kudla 
\paper The central critical value of a triple product $L$-function
\jour Annals of Math. 
\vol 133
\yr 1991
\pages 605--672
\endref

\ref\key{\harriskudlaIII} 
\bysame
\paper On a conjecture of Jacquet
\jour preprint, 2001, arXiv:math.NT/0111238
\endref



\ref\key{\hirzebruchzagier}
\by F. Hirzebruch and D. Zagier
\paper Intersection numbers of curves on Hilbert modular surfaces and modular forms
of Nebentypus
\jour Invent. Math.
\yr 1976
\vol 36
\pages 57--113
\endref

\ref\key{\howe}
\by R. Howe
\paper $\theta$--series and invariant theory
\inbook Proc. Symp. Pure Math.
\vol 33
\yr 1979
\pages 275--285
\endref

\ref\key{\howeps}
\by R. Howe and I. I. Piatetski-Shapiro
\paper Some examples of automorphic forms on $\roman{Sp}_4$
\jour Duke Math. J. 
\vol 50
\yr 1983
\pages 55--106
\endref

\ref\key{\hriljac}
\by P. Hriljac
\paper Heights and arithmetic intersection theory 
\jour Amer. J. Math.
\yr 1985
\vol 107
\pages 23--38
\endref

\ref\key{\jacquet}
\by H. Jacquet
\paper On the nonvanishing of some L-functions
\jour Proc. Indian Acad. Sci.
\vol 97
\yr 1987
\pages 117--155
\endref

\ref\key{\jacquetwald}
\bysame
\paper Sur un result de Waldspurger
\jour Ann. Scinet. \'Ec. Norm. Sup.
\vol 19
\yr 1986
\pages 185--229
\endref

\ref\key{\kitaoka}
\by Y. Kitaoka
\paper A note on local densities of quadratic forms
\jour Nagoya Math. J. 
\vol 92
\yr 1983
\pages 145--152
\endref 


\ref\key{\kottwitz} 
\by R. Kottwitz 
\paper Points on some Shimura varieties over finite fields 
\jour J.\ AMS 
{\bf 5} 
\yr 1992 
\pages 373--444 
\endref 


\ref\key{\annals}
\by S. Kudla
\paper Central derivatives of Eisenstein series and height pairings
\jour  Ann. of Math. 
\vol 146
\yr 1997 
\pages 545-646
\endref



\ref\key{\bourbaki}
\bysame
\paper Derivatives of Eisenstein series and generating functions for arithmetic cycles
\inbook S\'em. Bourbaki n${}^o$ 876  
\bookinfo Ast\'erisque 
\vol 276
\yr 2002
\pages 341--368
\endref

\ref\key{\Bints}
\bysame
\paper Integrals of Borcherds forms
\jour to appear in Compositio Math
\yr 
\endref

\ref\key{\kudlaICM} 
\bysame 
\paper Eisenstein series and arithmetic geometry
\jour Proc. of the ICM, Beijing, August 2002
\endref

\ref\key{\kudlaMSRI}
\bysame
\paper Special cycles and derivatives of Eisenstein series
\jour Proc. of Conference on Special Values of Rankin L-Series, MSRI, 
Dec. 2001, to appear.
\endref

\ref\key{\ariththeta}
\bysame
\paper An arithmetic theta function
\jour in preparation
\endref


\ref\key{\kmI} 
\by S. Kudla and J. Millson 
\paper The theta correspondence and harmonic forms I
\jour Math. Annalen 
\vol 274
\yr 1986
\pages 353--378
\endref

\ref\key{\krannals}
\by S. Kudla and S. Rallis
\paper A regularized Siegel-Weil formula: the first term identity  
\jour Annals of Math. 
\vol 140
\yr 1994
\pages 1--80
\endref

\ref\key{\krHB}
\by S. Kudla and M. Rapoport 
\paper Arithmetic Hirzebruch--Zagier cycles
\jour J. reine angew. Math. 
\vol 515 
\yr 1999
\pages 155--244
\endref 

\ref\key{\krinvent}
\bysame 
\paper Height pairings on Shimura curves and $p$-adic uniformization
\jour Invent. math.
\yr 2000
\vol 142
\pages 153--223
\endref

\ref\key{\krsiegel}
\bysame
\paper  Cycles on Siegel threefolds and derivatives of Eisenstein series
\jour Ann. Scient. \'Ec. Norm. Sup.
\vol 33
\yr 2000
\pages 695--756
\endref


\ref\key{\tiny}
\by S. Kudla, M. Rapoport and T. Yang 
\paper  On the derivative of an Eisenstein series of weight 1 
\jour Int. Math. Res. Notices, No.7 
\yr 1999 
\pages 347--385
\endref

\ref\key{\kryII}
\bysame
\paper Derivatives of Eisenstein series and Faltings heights
\jour preprint
\yr 2001
\vol
\pages
\endref

\ref\key{\kryIII}
\bysame
\paper Derivatives of Eisenstein series and Faltings heights II:  
modular curves
\jour in preparation
\endref

\ref\key{\kryIV}
\bysame
\paper On a Siegel--Eisenstein series
\jour in preparation
\endref

\ref\key{\doubling}
\bysame
\paper An arithmetic inner product formula
\jour in preparation
\endref



\ref\key{\kuehncrelle}
\by U. K\"uhn
\paper Generalized arithmetic intersection numbers
\jour J. reine angew. Math.
\yr 2001
\vol 534
\pages 209--236
\endref

\ref\key{\li}
\by J.-S. Li
\paper Non-vanishing theorems for the cohomology of certain 
arithmetic quotients
\jour J. reine angew. Math.
\vol428
\yr 1992
\pages 177--217
\endref

\ref\key{\maillotroessler}
\by V. Maillot and D. Roessler
\paper  Conjectures sur les d\'eriv\'ees logarithmiques des fonctions L 
d'Artin aux entiers n\'egatifs
\jour preprint
\yr 2001
\endref


\ref\key{\mcgraw}
\by W. J. McGraw
\paper On the rationality of vector-valued modular forms
\jour Math. Ann. 
\yr 2003
\pages DOI: 10.1007/s00208-003-0413-1
\endref

\ref\key{\niwa}
\by S. Niwa 
\paper Modular forms of half integral weight  and the integral of certain theta functions
\vol 56
\jour Nagoya Math. J. 
\yr 1975
\pages 147--161
\endref



\ref\key{\odatsuzuki}
\by T. Oda and M. Tsuzuki
\paper Automorphic Green functions associated to secondary 
spherical harmonics
\jour preprint
\yr 2000
\endref

\ref\key{\pswald}
\by I. I. Piatetski-Shapiro
\paper Work of Waldspurger
\inbook Lie Group Representations II
\bookinfo Lecture Notes in Math. 1041
\yr 1984
\pages 280--302
\endref 


\ref\key{\psrallis}
\by I. I. Piatetski-Shapiro and  S. Rallis  
\book  $L\!$-functions for classical groups
\bookinfo Lecture Notes in Math.
\publ Springer-Verlag
\publaddr New York
\vol  1254 
\yr  1987 
\pages  1--52 
\endref

\ref\key{\rallisinnerprod}
\by S. Rallis 
\paper  Injectivity properties of liftings associated to Weil representations
\jour  Compositio Math. 
\vol  52
\yr  1984
\pages  139--169
\endref

\ref\key{\saito}
\by H. Saito
\paper On Tunnell's formula for characters of $\GL_2$
\jour Compositio Math
\vol 85
\yr 1993
\pages 99--108
\endref

\ref\key{\siegel}
\by C. L. Siegel
\book Lectures on quadratic forms
\publ Tata Institute 
\publaddr Bombay
\yr 1957
\endref

\ref\key{\shimurahalf}
\by G. Shimura
\paper On modular forms of half integral weight
\jour Ann. of Math.   
\vol 97  
\yr 1973
\pages 440--481
\endref

\ref\key{\shimuraconf}
\bysame
\paper Confluent hypergeometric functions on tube domains
\jour Math. Ann. 
\vol 260
\yr 1982 
\pages 269-302
\endref

\ref\key{\soulebook} 
\by C. Soul\'e, D. Abramovich, J.-F. Burnol, and J.Kramer
\book Lectures on Arakelov Geometry, 
\bookinfo Cambridge Stud. Adv. Math., vol 33 
\publ Cambridge U. Press.
\yr 1992.
\endref

\ref\key{\tunnell}
\by J. Tunnell
\paper Local $\e$--factors and characters of $\GL_2$
\jour Amer. J. of Math.
\vol 105
\yr 1983
\pages 1277--1308
\endref


\ref\key{\waldshimura}
\by J.-L. Waldspurger
\paper Correspondance de Shimura
\jour J. Math. Pures Appl.
\yr 1980
\vol 59
\pages 1--132
\endref

\ref\key{\waldfourier}
\bysame
\paper Sur les coefficients de Fourier des formes modulaires
de poids demi-entier
\jour J. Math. Pures Appl.
\vol 60
\pages 375--484
\yr 1981
\endref


\ref\key{\waldsurvey}
\bysame
\paper Correspondance de Shimura
\inbook S\'eminare de Th\'eorie des Nombres, Paris 1979--80
\bookinfo Progr. Math., 12
\publ Birkh\"auser
\publaddr Boston
\yr 1981
\pages 357--369
\endref

\ref\key{\waldcentral}
\bysame
\paper Sur les valeurs de certaines fonctions L automorphes en leur 
centre de sym\'etrie
\jour Compositio Math
\yr 1985
\vol 54
\pages 173--242
\endref

\ref\key{\waldspurgerLHD}
\bysame
\paper Demonstration d'une conjecture de duality de Howe dans le case 
p-adique, $p\ne2$
\inbook Festschrift in honor of Piatetski-Shapiro, vol 2
\bookinfo Israel Math. Conf. Proc. 
\yr 1990
\pages 267--234
\endref


\ref\key{\waldquaternion}
\bysame
\paper Correspondances de Shimura et quaternions
\jour Forum Math. 
\vol 3
\yr 1991
\pages 219--307
\endref

\ref\key{\yangden} 
\by T. Yang
\paper\rm  An explicit formula for local densities of quadratic forms
\jour\it J. Number Theory
\vol 72
\yr 1998
\pages 309--356
\endref


\ref\key{\yangiccm}
\bysame
\paper The second term of an Eisenstein series
\jour Proc. of the ICCM, (to appear)
\endref

\ref\key{\yangMSRI}
\bysame
\paper Faltings heights and the derivative of Zagier's Eisenstein series
\jour Proc. of MSRI workshop on Heegner points, preprint (2002)
\endref

\ref\key{\zagierII}
\by D. Zagier
\paper\rm  Nombres de classes et formes modulaires de poids 3/2
\jour\it C. R. Acad. Sc. Paris
\yr 1975
\vol 281
\pages 883--886
\endref


\ref\key{\zagier}
\bysame
\paper\rm  Modular points, modular curves, modular surfaces and modular forms
\inbook Lecture Notes in Math. 1111
\yr 1985
\pages 225--248
\publ Springer
\publaddr Berlin
\endref

\ref\key{\zhang}
\by Shou-Wu Zhang
\paper Gross--Zagier formula for $\roman{GL}_2$
\jour Asian J. of Math. 
\vol 5
\yr 2001
\pages 183--290
\endref
\bye